\documentclass[reqno]{amsart}
\usepackage{setspace}
\usepackage{graphicx}
\usepackage{amsmath}
\usepackage{amsfonts}
\usepackage{amssymb}
\twocolumn \textwidth=16.5cm\textheight=24.6cm
\begin{document}

\noindent
{\LARGE \bf Can You Pave the Plane Nicely}

\medskip\noindent
{\LARGE\bf with Identical Tiles}

\bigskip\medskip
\noindent
{\large\it In Memory of Professor Peter M. Gruber}

\bigskip\medskip
\hfill{\bf Chuanming Zong}

\vspace{1cm}
\noindent
{\bf Abstract.} Everybody knows that identical regular triangles or squares can tile the whole plane. Many people know that
identical regular hexagons can tile the plane properly as well. In fact, even the bees know and use this fact! Is there any other convex
domain which can tile the Euclidean plane? Of course, there is a long list of them! To find the list and to show the completeness of the
list is a unique drama in mathematics, which has lasted for more than one century and the completeness of the list has
been mistakenly announced not only once! Up to now, the list consists of triangles, quadrilaterals, fifteen types of pentagons, and three types of hexagons. In 2017, Micha\"el Rao announced a computer proof for the completeness of the list. Meanwhile,
Qi Yang and Chuanming Zong made a series of unexpected discoveries in multiple tilings in the Euclidean plane. For examples, besides
parallelograms and centrally symmetric hexagons, there is no other convex domain which can form any two-, three- or four-fold translative tiling in the plane. However, there are two types of octagons and one type of decagons which can form nontrivial five-fold translative tilings. Furthermore, a convex domain can form a five-fold translative tiling of the plane if and only if it can form a five-fold lattice tiling. This paper tells the stories of these discoveries.


\vspace{1cm}
\noindent
{\Large\bf 1. Introduction}

\bigskip\noindent
Planar tilings is an ancient subject in our civilization. It has been considered in the arts by craftsmen since antiquity.
According to Gardner \cite{gard75}, the ancient Greeks knew that, among the regular polygons, only the triangle, the square
and the hexagon can tile the plane. Aristotle should know this fact since he made a similar claim in the space: {\it Among the five
Platonic solids, only the tetrahedron and the cube can tile the space}. Unfortunately, he made a mistake: {\it Identical regular
tetrahedra can not tile the whole space}!

The first recorded scientific investigation in tilings was made by Kepler. Assume that $\mathcal{T}$ is a tiling of the
Euclidean plane $\mathbb{E}^2$ by regular polygons. If the polygons are identical, the answer was already known to the ancient Greeks.
When different polygons are allowed, the situation becomes more complicated and more interesting. In particular, an edge-to-edge tiling
$\mathcal{T}$ by regular polygons is said to be of type $(n_1, n_2, \ldots , n_r)$ if each vertex ${\bf v}$ of $\mathcal{T}$ is surrounded
by a $n_1$-gon, a $n_2$-gon and so on in a cyclic order. Usually, they are known as {\it Archimedean tilings}. In 1619, Kepler \cite{kepl19}
enumerated all such tilings as $(3,3,3,3,3,$ $3)$, $(3,3,3,3,6)$, $(3,3,3,4,4)$, $(3,3,4,3,4)$, $(3,4,6,$ $4)$, $(3,6,3,6)$, $(3,12,12)$, $(4,4,4,4)$, $(4,6,12)$, $(4,8,8)$ and $(6,6,6)$. Beautiful illustrations of the Archimedean tilings can be found in many references such as \cite{grsh86, schu93}.

If ${\bf a}_1$, ${\bf a}_2$, $\ldots $, ${\bf a}_n$ are $n$ linearly independent vectors in $\mathbb{E}^n$, then the set
$$\Lambda =\left\{\sum z_i{\bf a}_i:\ z_i\in \mathbb{Z}\right\}$$
is an {\it $n$-dimensional lattice}. Clearly, lattices are the most natural periodic discrete sets in the plane and spaces. Therefore, many pioneering scientists like Kepler, Huygens, Ha\"uy and Seeber took lattice packings and lattice tilings as the scientific foundation for crystals. In 1885, the famous crystallographer Fedorov \cite{fedo85} discovered that: {\it A convex domain can form a lattice tiling of $\mathbb{E}^2$ if and only if it is a parallelogram or a centrally symmetric hexagon; a convex body can form a lattice tiling in $\mathbb{E}^3$ if and only if it is a parallelotope, a hexagonal prism, a rhombic dodecahedron, an elongated octahedron or a truncated octahedron.}

Usually, tilings allow very general settings, without restriction on the shapes of the tiles and the number of the different shapes. However,
to avoid complexity and confusion, in this paper we only deal with the tilings by identical convex polygon tiles. In particular, we call it a {\it translative tiling} if all the tiles are translates of each others, and call it a {\it lattice tiling} if it is a translative tiling and all the translative vectors together is a lattice.

In 1900, Hilbert \cite{hilb00} proposed a list of mathematical problems in his ICM lecture in Paris. As a generalized inverse of Fedorov's discovery, he wrote in his 18th problem that: {\it A fundamental region of each group of motions, together with the congruent regions arising from the group, evidently fills up space completely. The question arises: whether polyhedra also exist which do not appear as fundamental regions of groups of motions, by means of which nevertheless by a suitable juxtaposition of congruent copies a complete filling up of all space is possible.}

Hilbert proposed his problem in the space, perhaps he believed that there is no such domain in the plane. When Reinhardt started his Doctoral thesis at Frankfurt am Main in 1910s, Bieberbach suggested him to determine all the convex domains which can tile the whole plane and so that to verify that Hilbert's problem indeed has positive answer in the plane. This is the origin of the natural problem that {\it to determine all the two-dimensional convex tiles}. In 1917 Reinhardt was an assistant of Hilbert at G\"ottingen. Therefore they should have good chances to discuss this problem.

It is worth to mention that Bieberbach himself solved the first part of Hilbert's 18th problem in 1911.

\vspace{0.6cm}
\noindent
{\Large\bf 2. Reinhardt's List}

\bigskip\noindent
In 1918, Reinhardt \cite{rein18} made his Doctoral degree under the supervision of Bieberbach at Frankfurt am Main with a thesis \lq\lq{\it \"Uber die Zerlegung der Ebene in Polygone}". This is the first approach to characterize all the convex domains which can tile the whole plane. First, he studied the tiling networks and obtained an expression for the mean value of the vertices over faces. As a corollary of the formula, he obtained the following
result.

\medskip\noindent
{\bf Theorem 1.} {\it A convex $m$-gon can tile the whole plane $\mathbb{E}^2$ only if}
$$m\le 6.$$

In fact, as he and many other authors pointed out (see \cite{gard75,kers68,kers69,mann18,rein18}), this theorem can be easily deduced by {\it Euler's formula}
$$v-e+f=2,$$
where $v$, $e$ and $f$ stand for the numbers of the vertices, the edges and the faces of a polygonal division of the plane.

Apparently, two identical triangles can make a parallelogram and two identical quadrilaterals can make a centrally symmetric hexagon (see Figure 1). Thus, by Fedorov's theorem, identical triangles or quadrilaterals can always tile the plane nicely. However, it is easy to see that identical regular pentagons or some particular hexagons can not tile the plane. Then, Bieberbach's problem can be reformulated as:

\medskip
{\it What kind of convex pentagons or hexagons can tile the plane}?

\begin{figure}[!ht]
\centering
\includegraphics[scale=0.53]{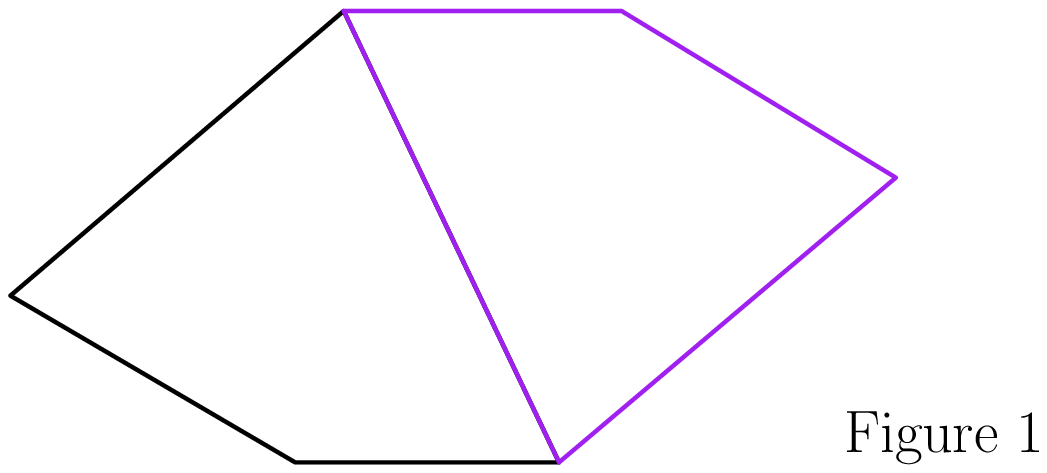}
\end{figure}

Let $P_n$ denote a convex $n$-gon with vertices ${\bf v}_1$, ${\bf v}_2$, $\ldots ,$ ${\bf v}_n$ in the anti-clock order, let $G_i$ denote the edge
with ends ${\bf v}_{i-1}$ and ${\bf v}_i$, where ${\bf v}_0={\bf v}_n$, let $\ell_i$ denote the length of $G_i$, and let $\alpha_i$ denote the inner angle of $P_n$ at ${\bf v}_i$.

Reinhardt's thesis obtained the following complete solution to the hexagon case of Bieberbach's problem.

\medskip\noindent
{\bf Theorem 2 (Reinhardt \cite{rein18}).} {\it A convex hexagon $P_6$ can tile the whole plane $\mathbb{E}^2$ if and only if it satisfies one of the three groups of conditions:

\noindent\begin{enumerate}
\item[\bf (1).] $\alpha_1+\alpha_2+\alpha_3=2\pi$ and $\ell_1=\ell_4$.
\item[\bf (2).] $\alpha_1+\alpha_2+\alpha_4=2\pi$, $\ell_1=\ell_4$, and $\ell_3=\ell_5$.
\item[\bf (3).] $\alpha_1=\alpha_3=\alpha_5={2\over 3}\pi$, $\ell_1=\ell_2$, $\ell_3=\ell_4$, and $\ell_5=\ell_6$.
\end{enumerate}}

\begin{figure}[!ht]
\centering
\includegraphics[scale=0.4]{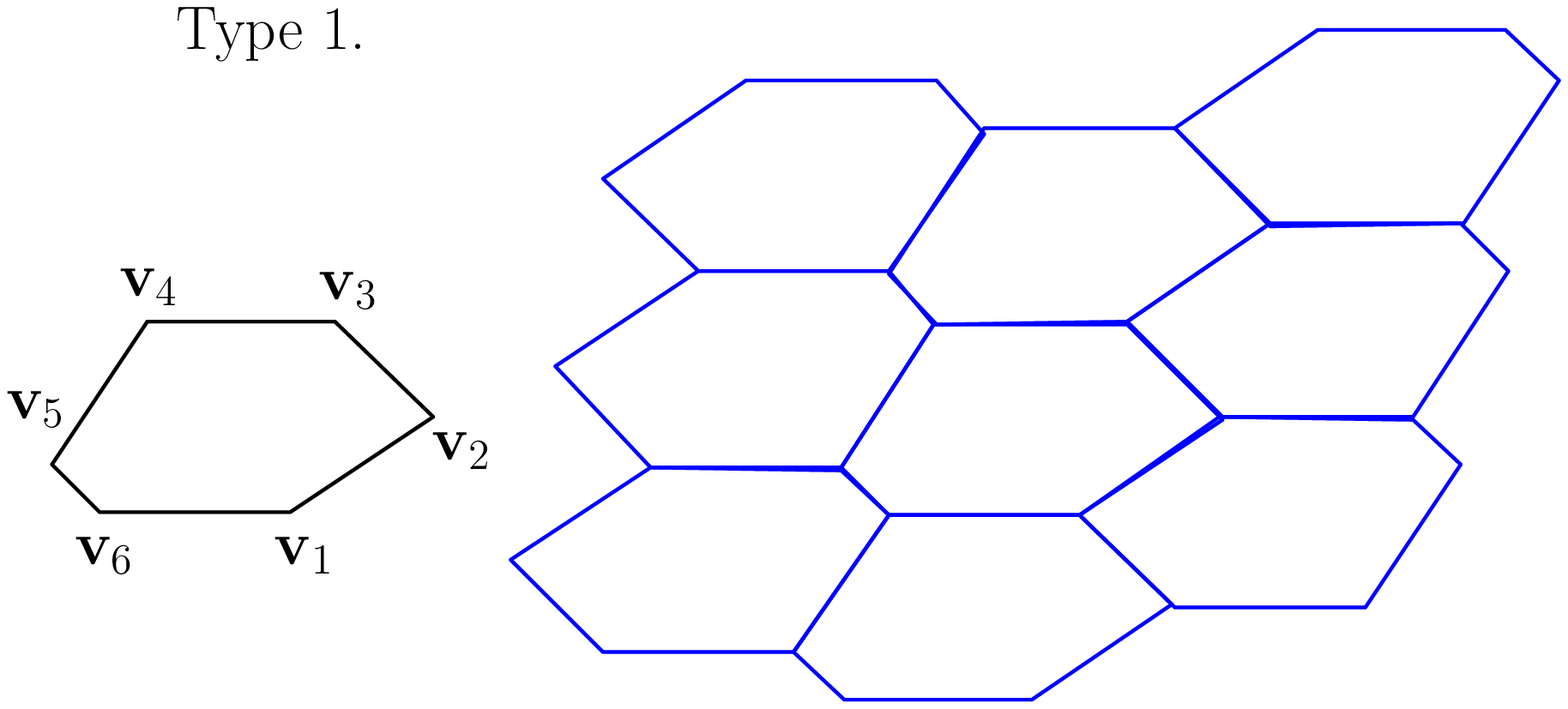}

\vspace{0.5cm}
\centering
\includegraphics[scale=0.37]{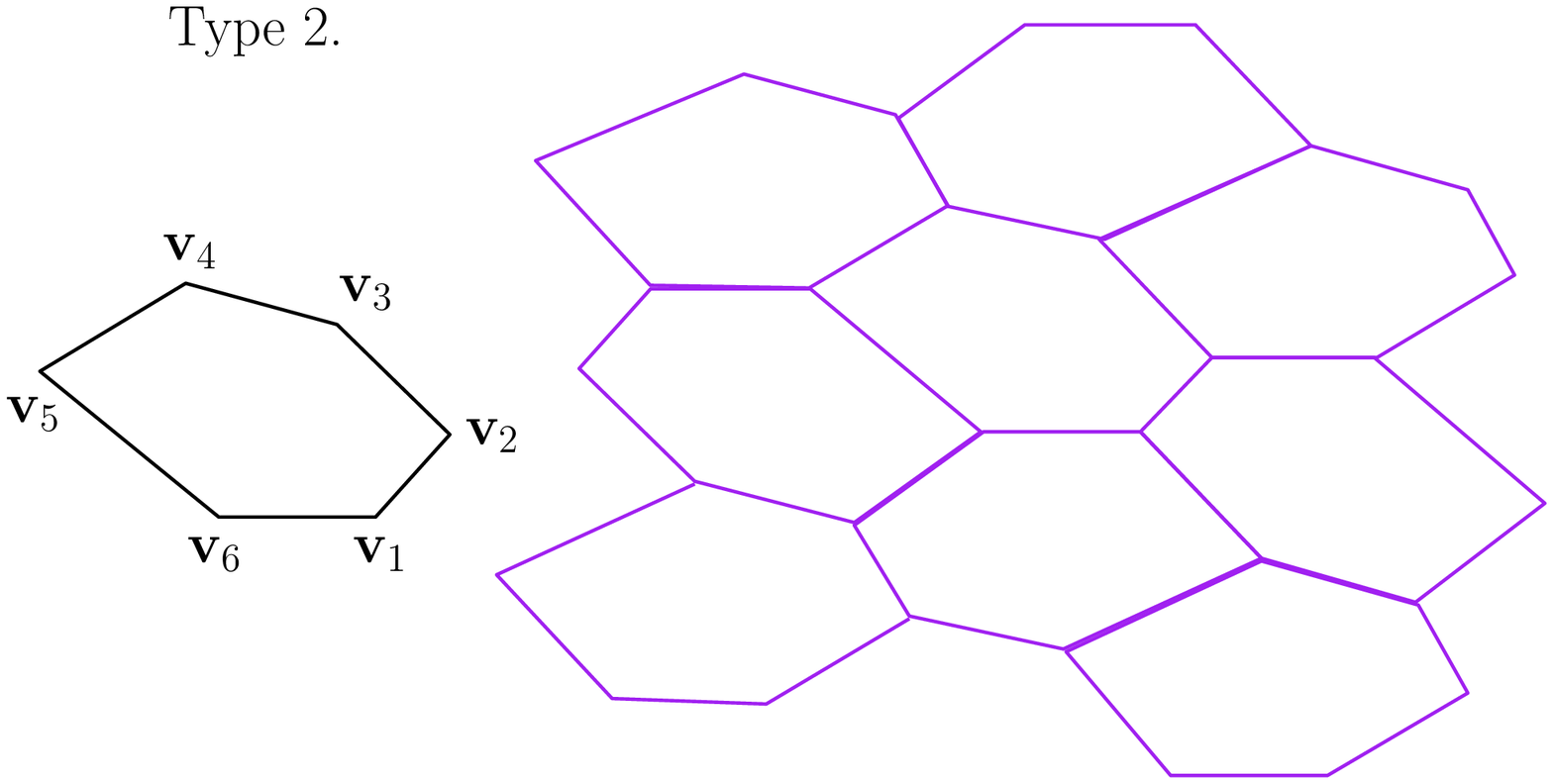}

\vspace{0.4cm}
\centering
\includegraphics[scale=0.40]{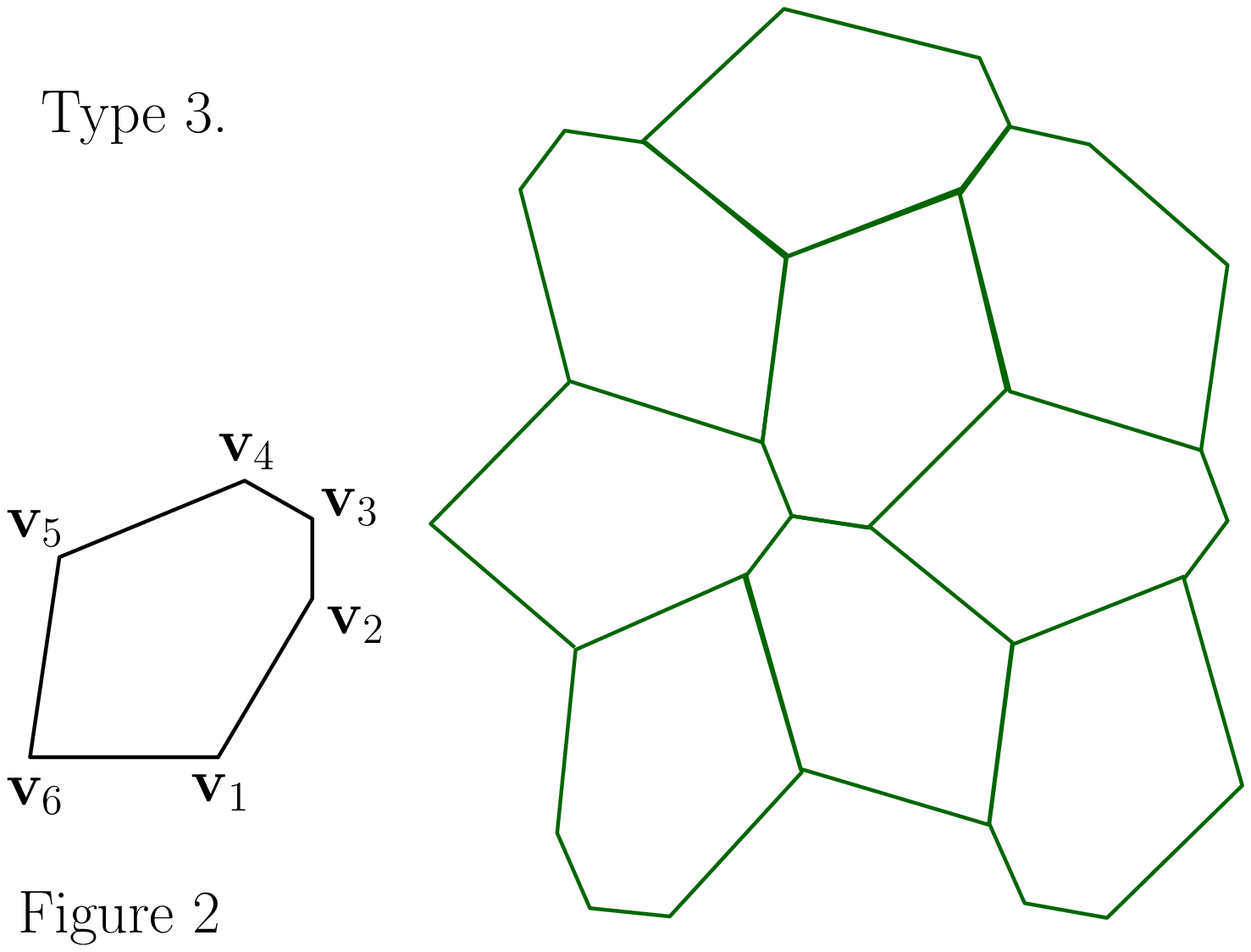}
\end{figure}

The if part of this theorem is relatively simple. It is illustrated by Figure 2. However, the only if part is much more complicated. Reinhardt deduced the only if part by considering six cases with respect to how many edges of the considered hexagon are equal. His proof was very sketchy and difficult to understand and check. It seems that he only considered the edge-to-edge tilings.

Fortunately, this theorem has been verified by several other authors. For example, without the knowledge of Reinhardt's thesis, in 1963 Bollob\'as made the following surprising observation, which guarantees the sufficiency of Reinhardt's consideration.

\medskip\noindent
{\bf Lemma 1 (Bollob\'as \cite{boll63}).} {\it If $\mathcal{T}$ is a tiling of the plane by identical convex hexagons and $\gamma $ is any given positive number, there is a square of edge length $\gamma $ in which the tiling is edge-to-edge and every vertex is surrounded by three hexagons.}

\medskip
For the pentagon tilings, by considering five cases with respect to how many edges are equal, Reinhardt obtained the following result.

\medskip\noindent
{\bf Theorem 3 (Reinhardt \cite{rein18}).} {\it A convex pentagon $P_5$ can tile the whole plane $\mathbb{E}^2$ if it satisfies one of the five groups of conditions:

\noindent\begin{enumerate}
\item[\bf (1).] $\alpha_1+\alpha_2+\alpha_3=2\pi$.
\item[\bf (2).] $\alpha_1+\alpha_2+\alpha_4=2\pi$ and $\ell_1=\ell_4$.
\item[\bf (3).] $\alpha_1=\alpha_3=\alpha_4={2\over 3}\pi$, $\ell_1=\ell_2$ and $\ell_4=\ell_3+\ell_5$.
\item[\bf (4).] $\alpha_1=\alpha_3={1\over 2}\pi$, $\ell_1=\ell_2$ and $\ell_3=\ell_4$.
\item[\bf (5).] $\alpha_1={1\over 3}\pi$, $\alpha_3={2\over 3}\pi$, $\ell_1=\ell_2$ and $\ell_3=\ell_4$.
\end{enumerate}}

\begin{figure}[!ht]
\centering
\includegraphics[scale=0.4]{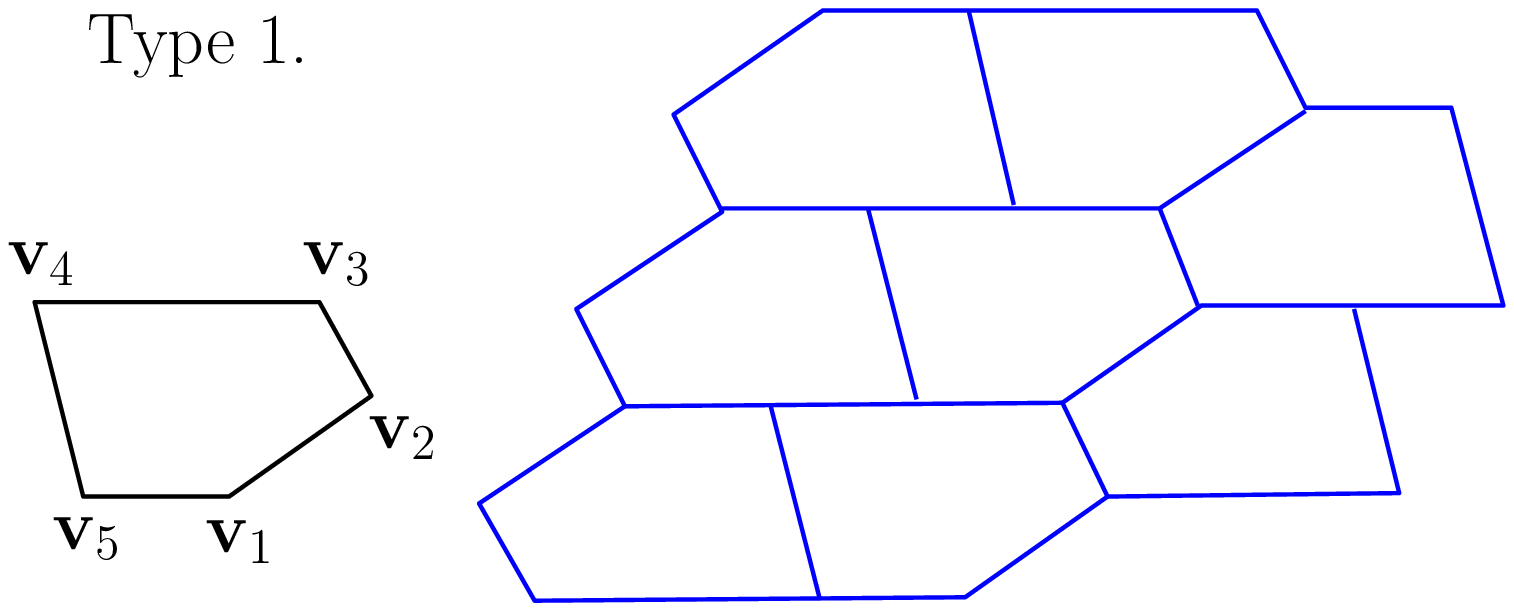}

\vspace{0.4cm}
\centering
\includegraphics[scale=0.4]{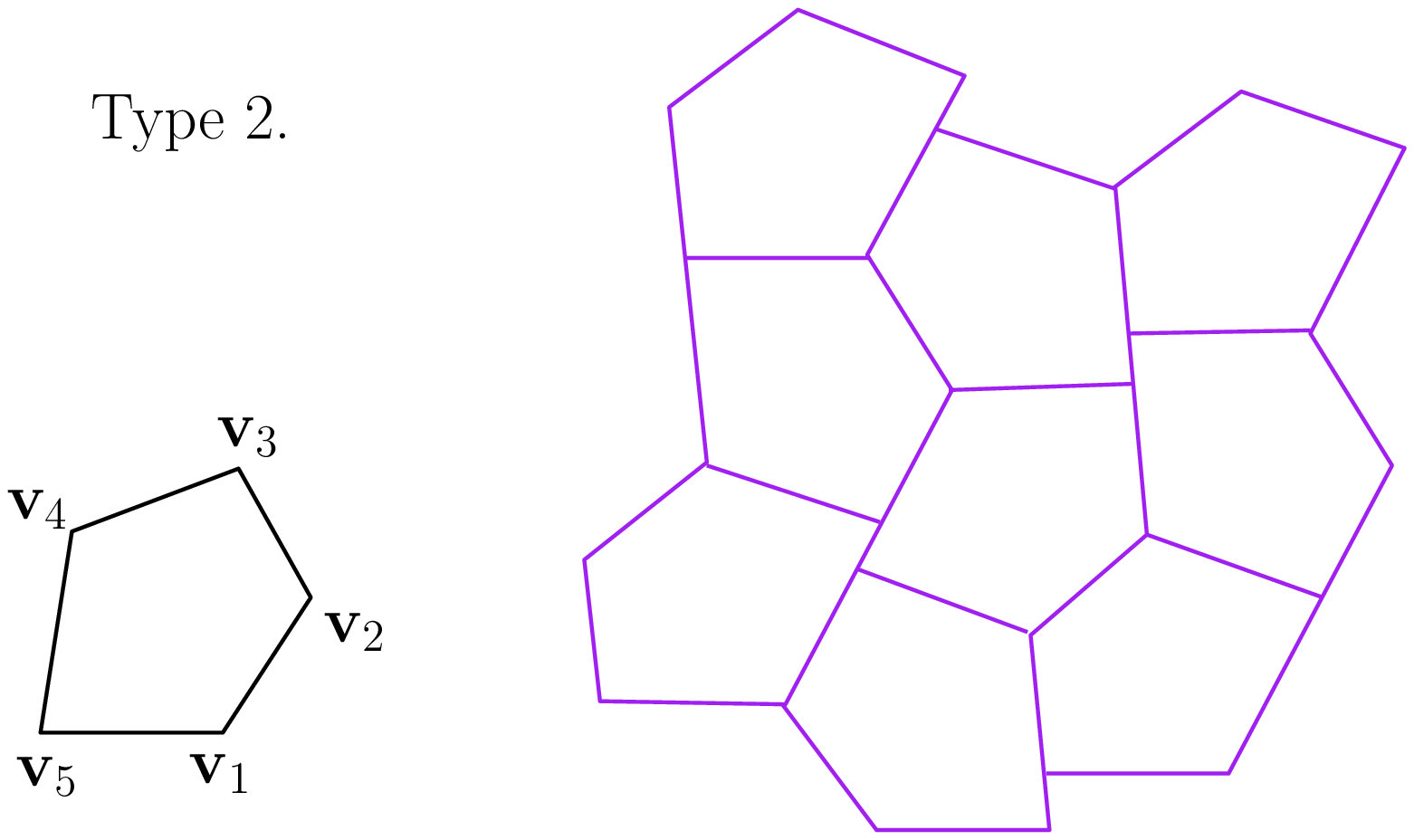}

\vspace{0.55cm}
\centering
\includegraphics[scale=0.4]{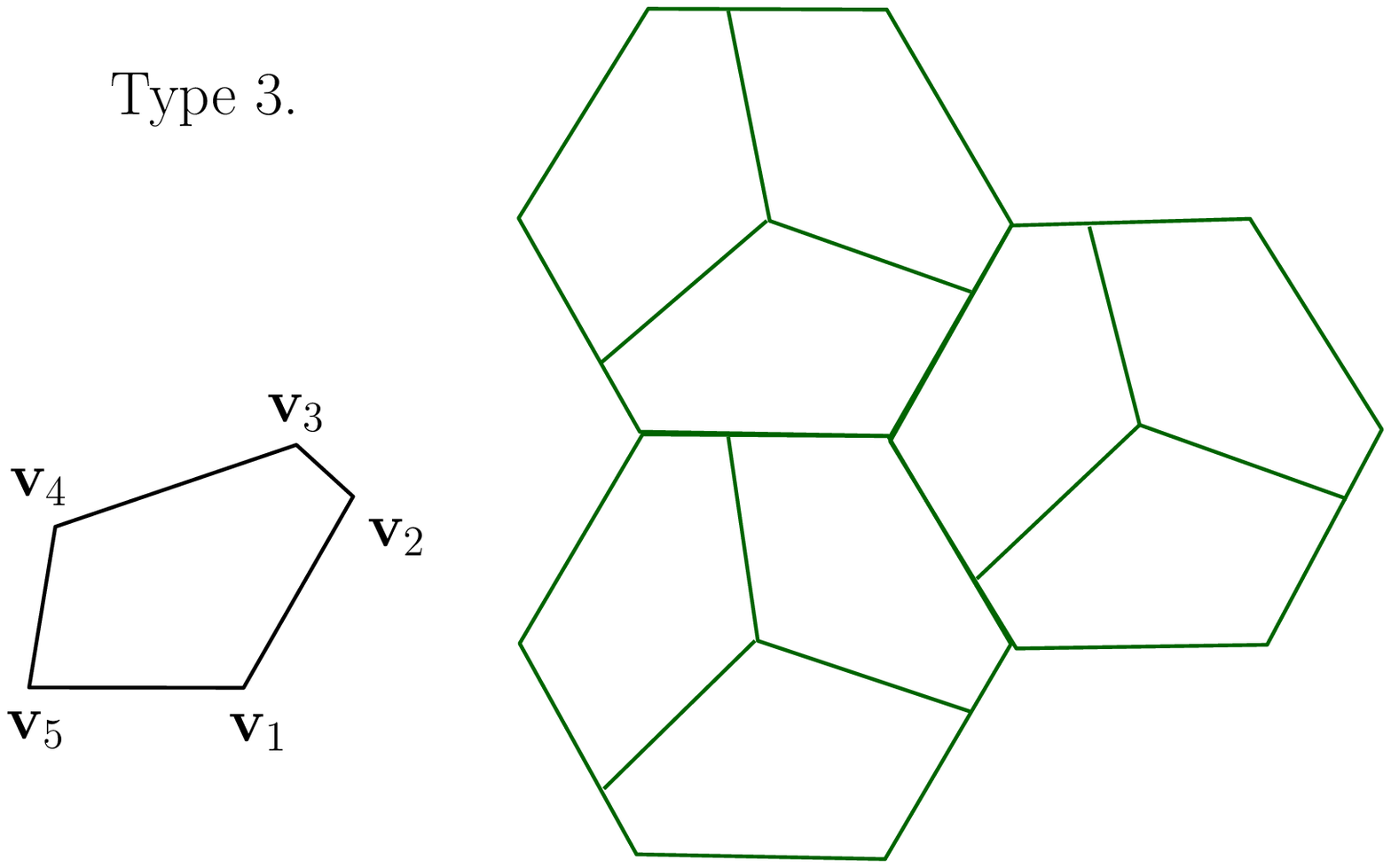}

\vspace{0.4cm}
\centering
\includegraphics[scale=0.4]{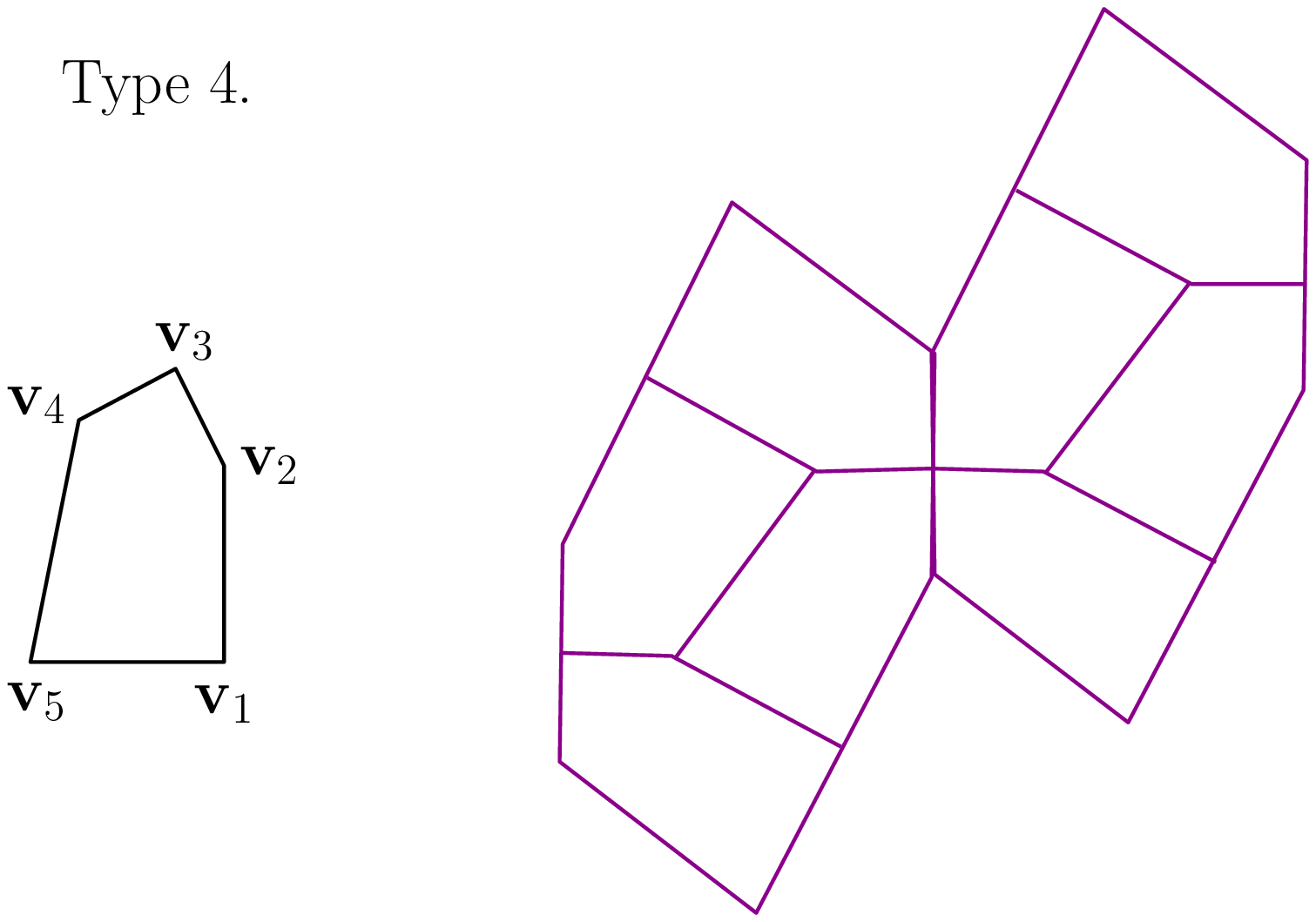}

\vspace{0.4cm}
\centering
\includegraphics[scale=0.4]{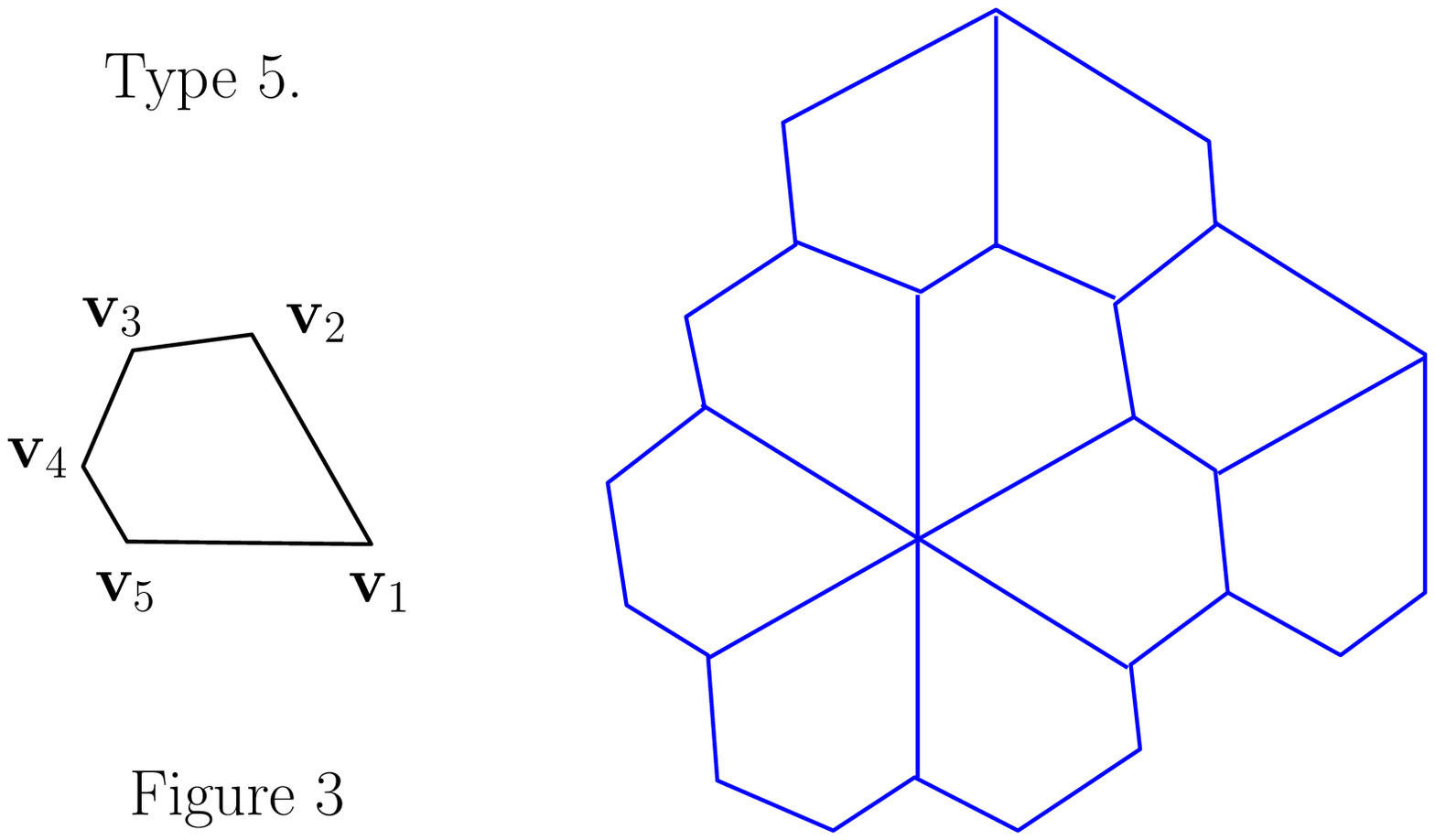}
\end{figure}

\smallskip
To figure out the list is nontrivial. However, as shown in Figure 3, it is easy to check that all the pentagons listed in Theorem 3 indeed can tile the plane. Reinhardt himself did not claim the completeness of the pentagon tile list. However, according to Gardner \cite{gard75} {\it it is quite clear that Reinhardt and everyone else in the field thought that the Reinhardt pentagon list was probably complete.}

As it was observed by Reinhardt \cite{rein18} that, all triangles, quadrilaterals, the three types of hexagons listed in Theorem 1 and the five classes of pentagons listed in Theorem 2 are indeed fundamental domains of some groups of motions. Both Hilbert and Bieberbach should be happy to know this.

Unfortunately, in 1928 Reinhardt \cite{rein28} discovered a (non-convex) three-dimensional polytope which can form a tiling in the space but is not the fundamental domain of any group of motions! This is the first counter-example to the second part of Hilbert's 18th problem.

Inspired by Reinhardt's discovery, in 1935 Heesch \cite{hees35} obtained a two-dimensional counter-example (see Figure 4) to Hilbert's problem.

\begin{figure}[!ht]
\centering
\includegraphics[scale=0.43]{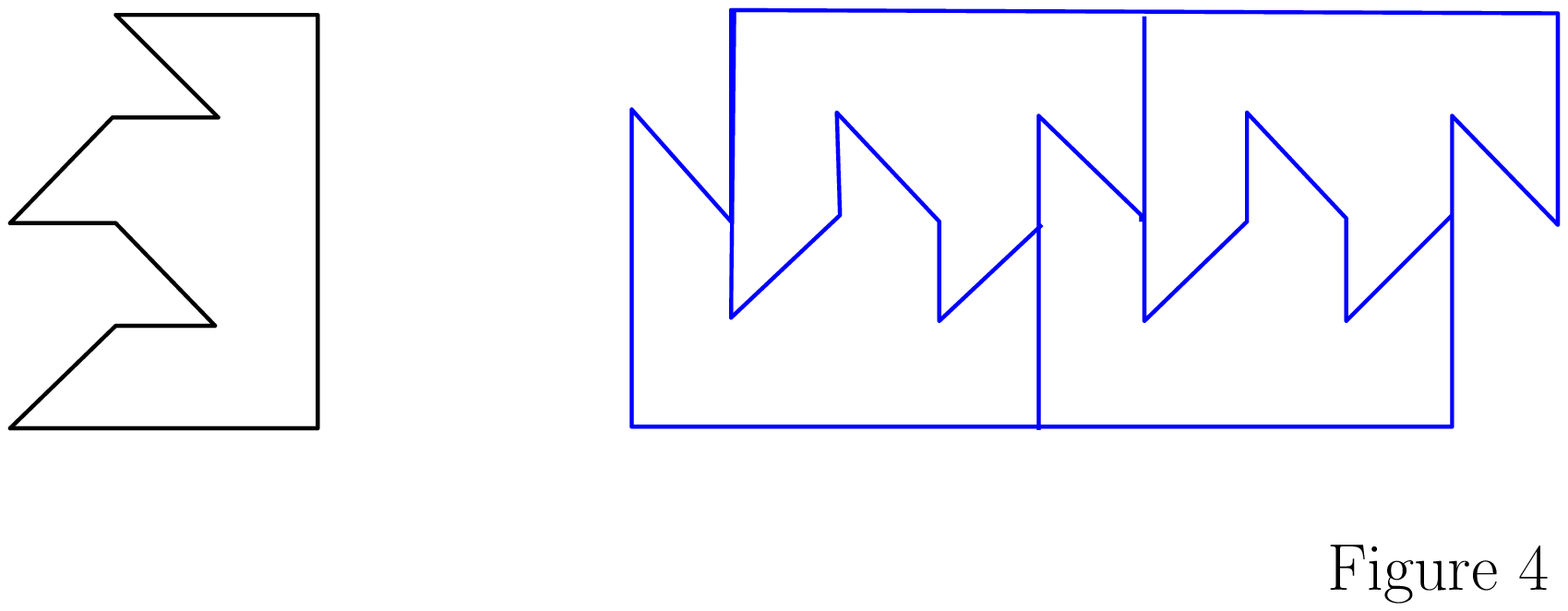}
\end{figure}

Thirty years later, Heesch and Kienzle \cite{heki63} presented a rather detailed treatment of plane tilings, including non-convex tiles. No new convex tile was discovered. It was claimed that the treatment given was complete.

\vspace{0.6cm}
\noindent
{\Large\bf 3. An End, Or A New Start}

\bigskip\noindent
In 1968, fifty years after Reinhardt's pioneering thesis, Kershner surprisingly discovered three new classes of pentagons which can pave the whole plane without gap and overlapping.

\medskip\noindent
{\bf Theorem 4 (Kershner \cite{kers68}).} {\it A convex pentagon $P_5$ can tile the whole plane $\mathbb{E}^2$ if it satisfies one of the three groups of conditions:

\noindent\begin{enumerate}
\item[\bf (6).] $\alpha_1+\alpha_2+\alpha_4=2\pi$, $\alpha_1=2\alpha_3$, $\ell_1=\ell_2=\ell_5$, and $\ell_3=\ell_4$.
\item[\bf (7).] $2\alpha_2+\alpha_3=2\alpha_4+\alpha_1=2\pi$ and $\ell_1=\ell_2=\ell_3=\ell_4$.
\item[\bf (8).] $2\alpha_1+\alpha_2=2\alpha_4+\alpha_3=2\pi$ and $\ell_1=\ell_2=\ell_3=\ell_4$.
\end{enumerate}}

\smallskip
According to Kershner \cite{kers69}, had been intrigued by this problem for some 35 years, he finally discovered a method of classifying the possibilities for pentagons in a more convenient way than Reinhardt's to yield an approach that was humanly possible to carry to completion. Unfortunately, neither \cite{kers68} nor \cite{kers69} contains any hint of his method. Of course, the three classes of new pentagon tiles were indeed surprising, though verifications are simple (see Figure 5).

In the introduction of \cite{kers69}, Kershner stated that: {\it The author has recently succeeded in carrying through a complete determination in the special case of convex pentagons. This paper contains a few critical proofs and a complete statement of the results, but not a complete proof for the excellent reason that a complete proof would require a rather large book.}

\begin{figure}[!ht]
\centering
\includegraphics[scale=0.42]{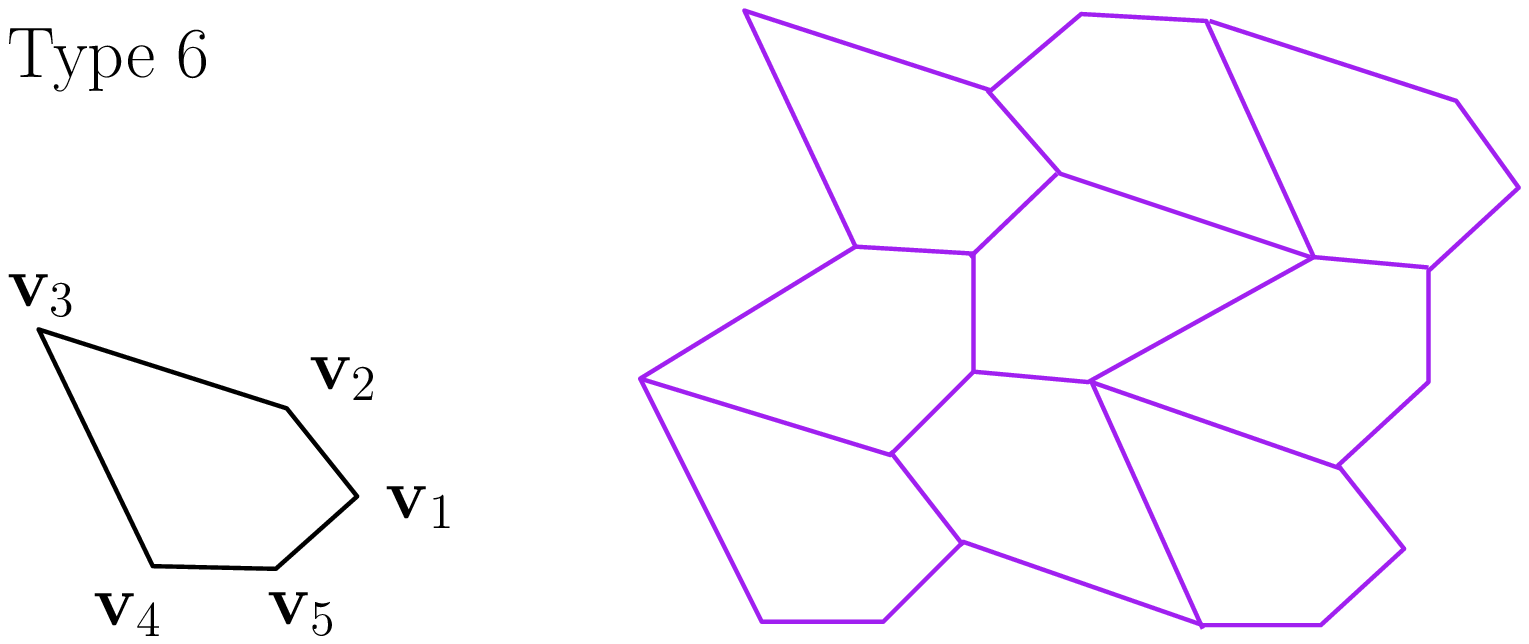}

\vspace{0.5cm}
\centering
\includegraphics[scale=0.4]{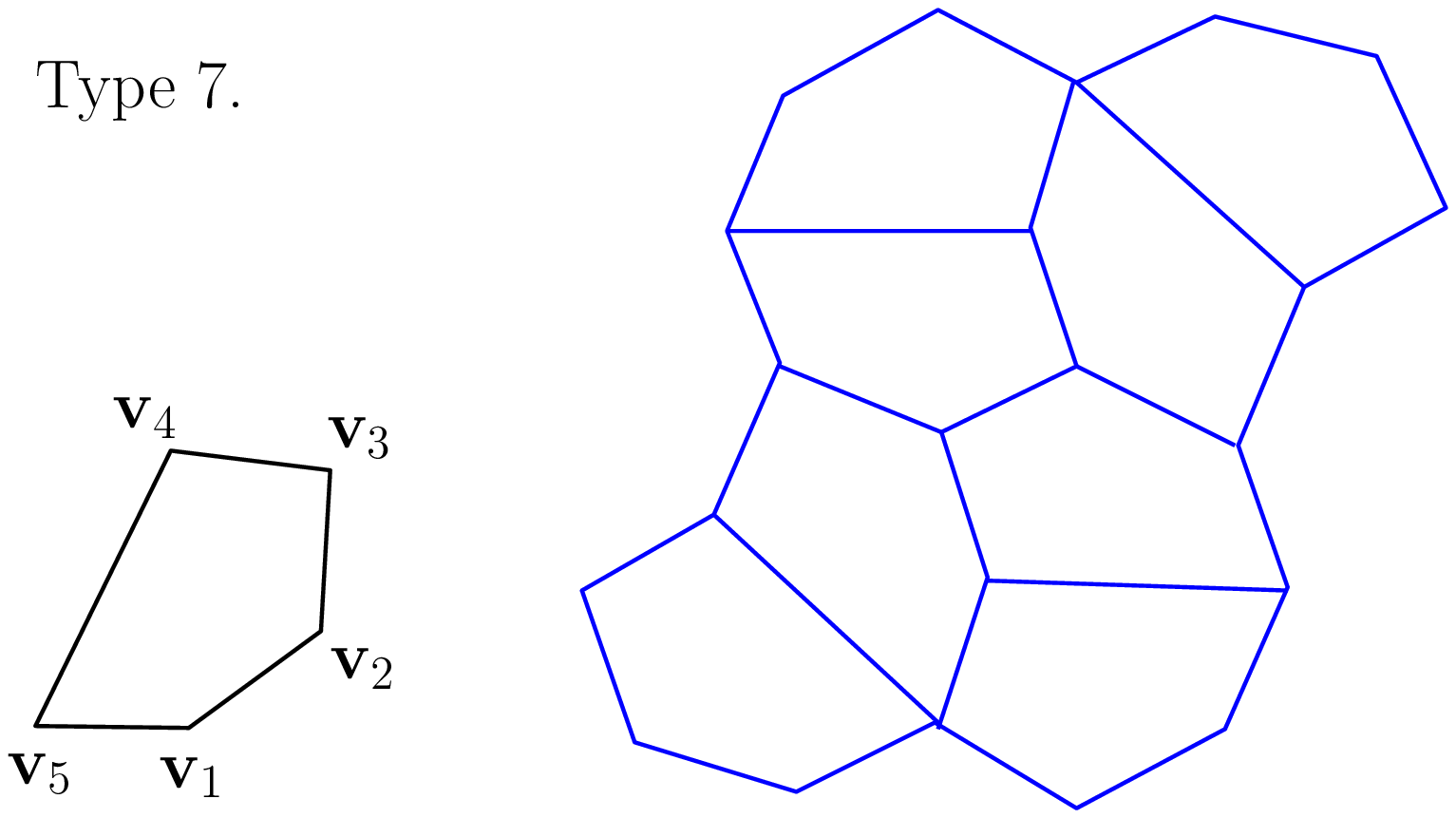}

\vspace{0.5cm}
\centering
\includegraphics[scale=0.4]{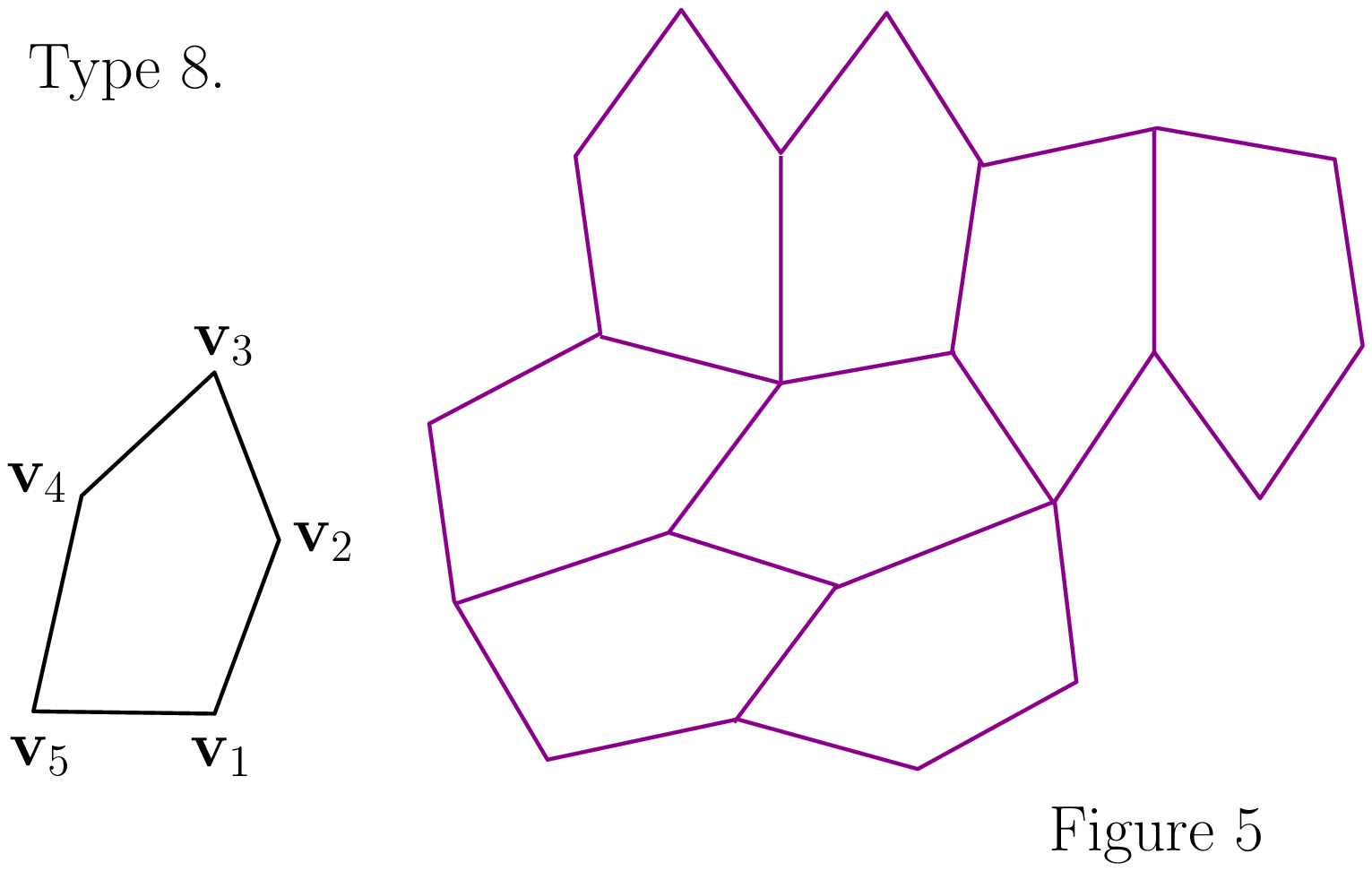}
\end{figure}

\smallskip\noindent
{\bf Remark 1.} The pentagons of Types 6-8 are counter-examples to the second part of Hilbert's 18th problem! In other words, they can tile the whole plane, nevertheless they are not the fundamental domains of any group of motions. Inductively, $n$-dimensional cylinder counter-examples can be constructed from $(n-1)$-dimensional ones. Needless to say, all Hilbert, Bieberbach, Reinhardt, Heesch and others should be surprised by Kershner's elegant examples! In fact, Kershner himself did not mention this fact in his papers. Perhaps he overlooked it. This fact has been mentioned in many books and survey papers, see \cite{bmp,egh,grsh80,schu93}.

\medskip
In 1975, Reinhardt and Kershner's discoveries were introduced by Martin Gardner, a famous scientific writer, at the mathematical games column of the Scientific American magazine. Since then the tiling problem has stimulated the interests of many amateurs. Surprisingly, they even have made remarkable contributions to this problem.

Soon after Gardner's popular paper, based on the known tiling pattern by octagons and squares together as shown by Figure 6, a computer scientist Richard Jammes III discovered a class of new tiles.

\begin{figure}[!ht]
\centering
\includegraphics[scale=0.42]{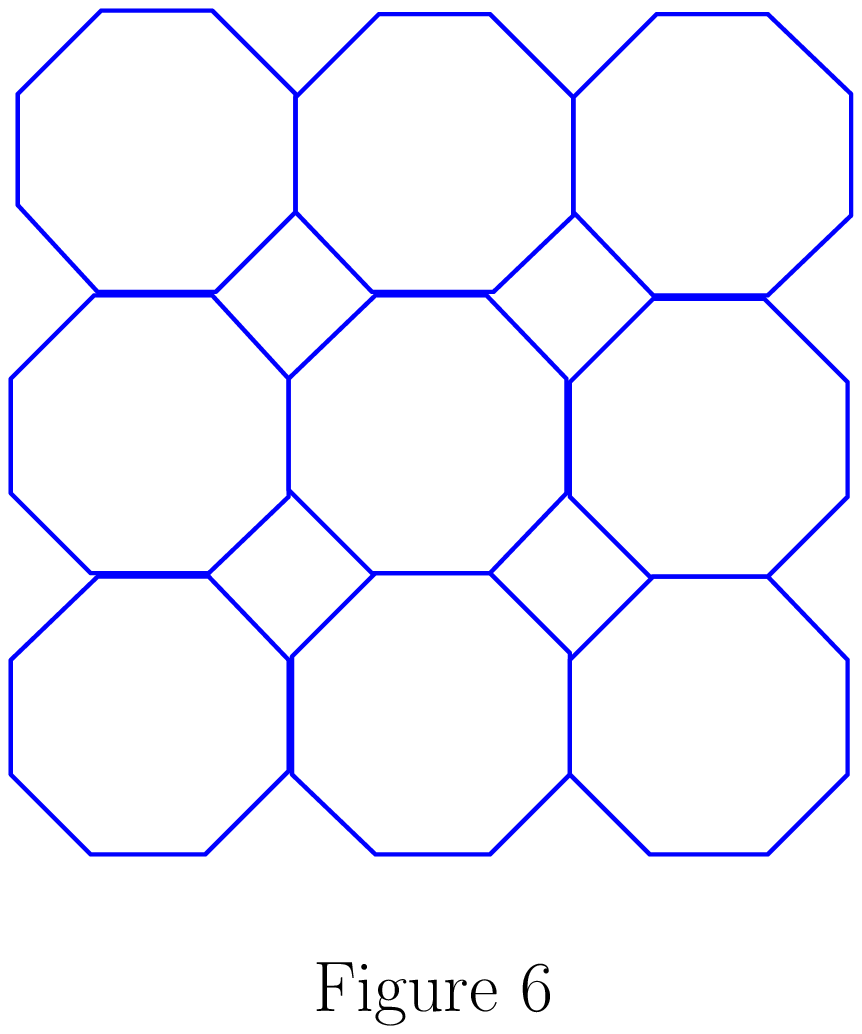}
\end{figure}

\medskip\noindent
{\bf Theorem 5 (James \cite{jame75}).} {\it A convex pentagon $P_5$ can tile the whole plane $\mathbb{E}^2$ if it satisfies the following group of conditions:

\noindent\begin{enumerate}
\item[\bf (9).] $\alpha_5={\pi \over 2}$, $\alpha_1+\alpha_4=\pi$, $2\alpha_2-\alpha_4=2\alpha_3+\alpha_4=\pi$ and $\ell_1=\ell_2+\ell_4=\ell_5$.
\end{enumerate}}

\smallskip
This result can be easily verified by argument based on Figure 7. In principle, Lemma 1 guarantees that every hexagon tiling is edge-to-edge.
However, James' discovery shows that this is no longer true in some pentagon tilings. Theorem 5 also served to point out that Kershner had taken edge-to-edge as a hidden assumption in his consideration.

\begin{figure}[!ht]
\centering
\includegraphics[scale=0.38]{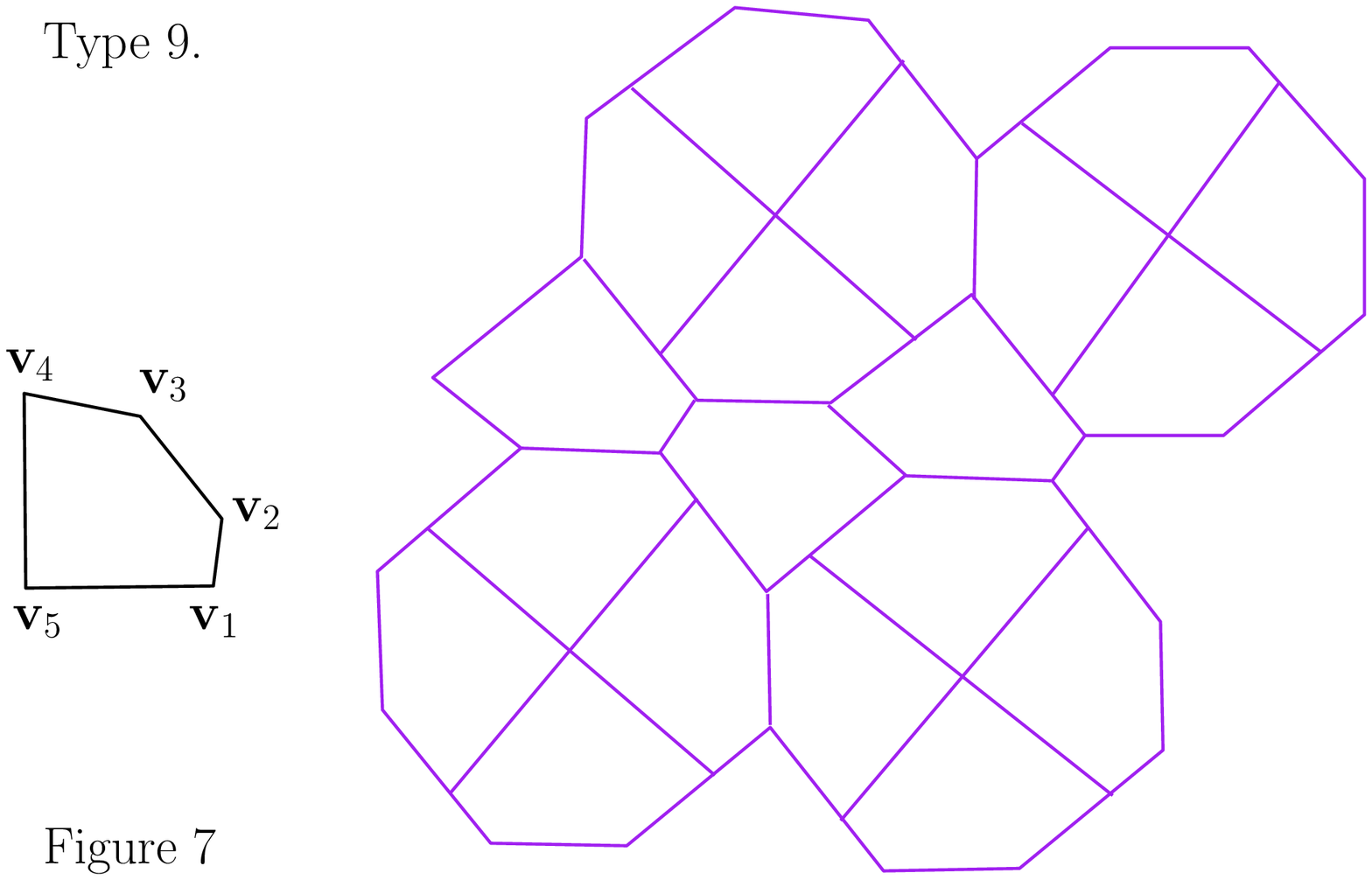}
\end{figure}

\smallskip
Meanwhile, Marjorie Rice made a true astonishing news. First, she was a true amateur. According to Schattschneider \cite{scha78,wolc}, Rice had no mathematical training beyond \lq\lq the bare minimum they required $\cdots$ in high school over 35 years ago". Second, even so she was able to consider the problem with a systematic method based on the possible local structures of the pentagon tilings at a given vertex. By dealing with
more than sixty cases, four types of new pentagon tiles were discovered!

\medskip\noindent
{\bf Theorem 6 (Rice \cite{scha78}).} {\it A convex pentagon $P_5$ can tile the whole plane $\mathbb{E}^2$ if it satisfies one of the four groups of conditions:

\noindent\begin{enumerate}
\item[\bf (10).] $\alpha_2+2\alpha_5=2\pi$, $\alpha_3+2\alpha_4=2\pi$ and $\ell_1=\ell_2=\ell_3=\ell_4$.
\item[\bf (11).] $\alpha_1={\pi \over 2}$, $\alpha_3+\alpha_5=\pi $, $2\alpha_2+\alpha_3=2\pi$ and $2\ell_1+\ell_3=\ell_4=\ell_5$.
\item[\bf (12).] $\alpha_1={\pi \over 2}$, $\alpha_3+\alpha_5=\pi $, $2\alpha_2+\alpha_3=2\pi$ and $2\ell_1=\ell_3+\ell_5=\ell_4$.
\item[\bf (13).] $\alpha_1=\alpha_3={\pi\over 2}$, $2\alpha_2+\alpha_4=2\alpha_5+\alpha_4=2\pi $, $\ell_3=\ell_4$ and $2\ell_3=\ell_5$.
\end{enumerate}}

\smallskip
It is routine to verify this theorem based on Figure 8. Nevertheless, it is rather surprising to notice that the tilings produced by the pentagons of type 10 are edge-to-edge, which was missed by both Reinhardt and Kershner. It is even more surprising that all the pentagons of types 9-13 are counter-examples to Hilbert's problem as well (see \cite{grsh80}). In other words, they can tile the whole plane, however they are not the fundamental domains of any group of motions.

\begin{figure}[!ht]
\centering
\includegraphics[scale=0.42]{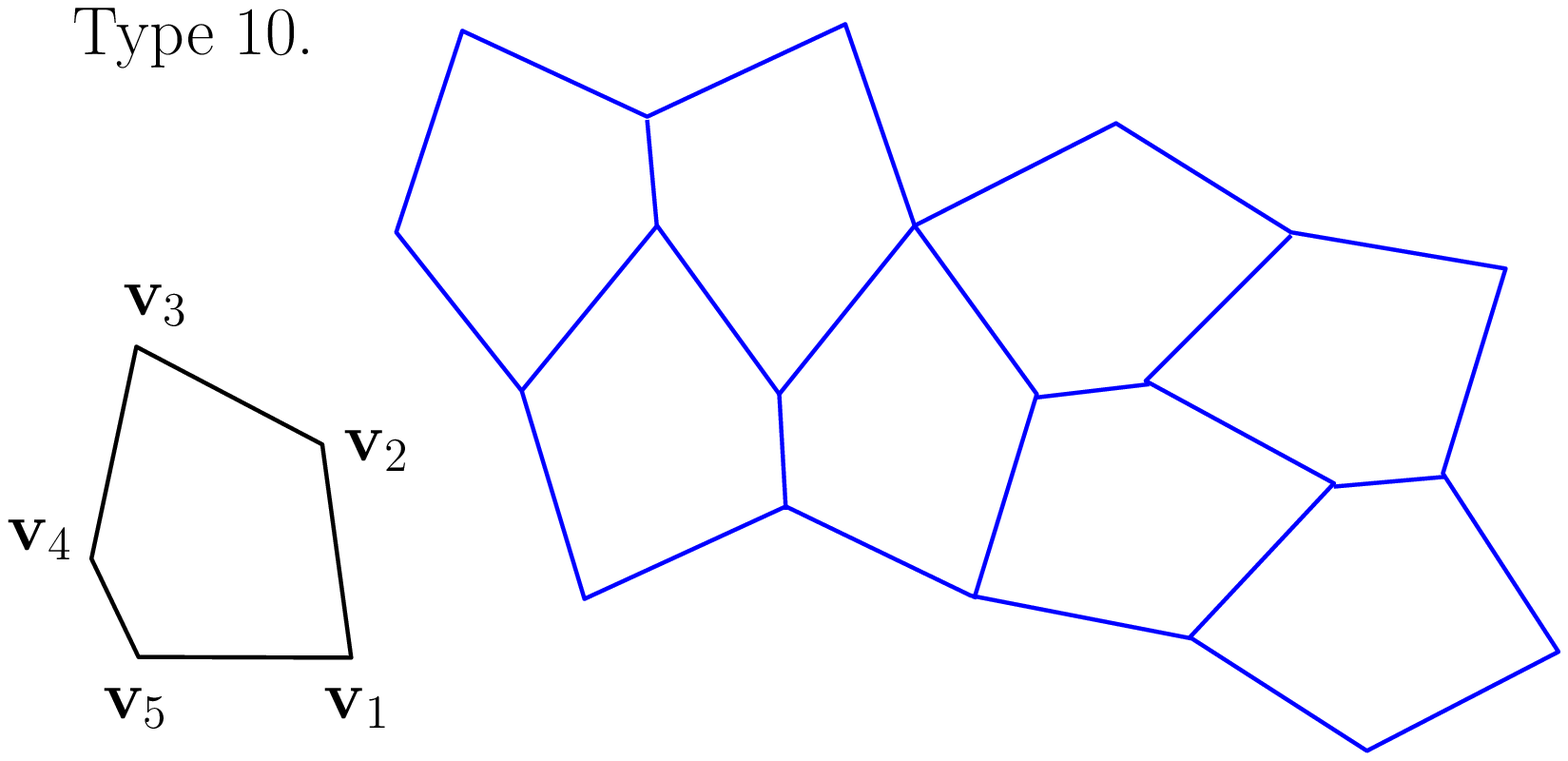}

\vspace{0.45cm}
\centering
\includegraphics[scale=0.45]{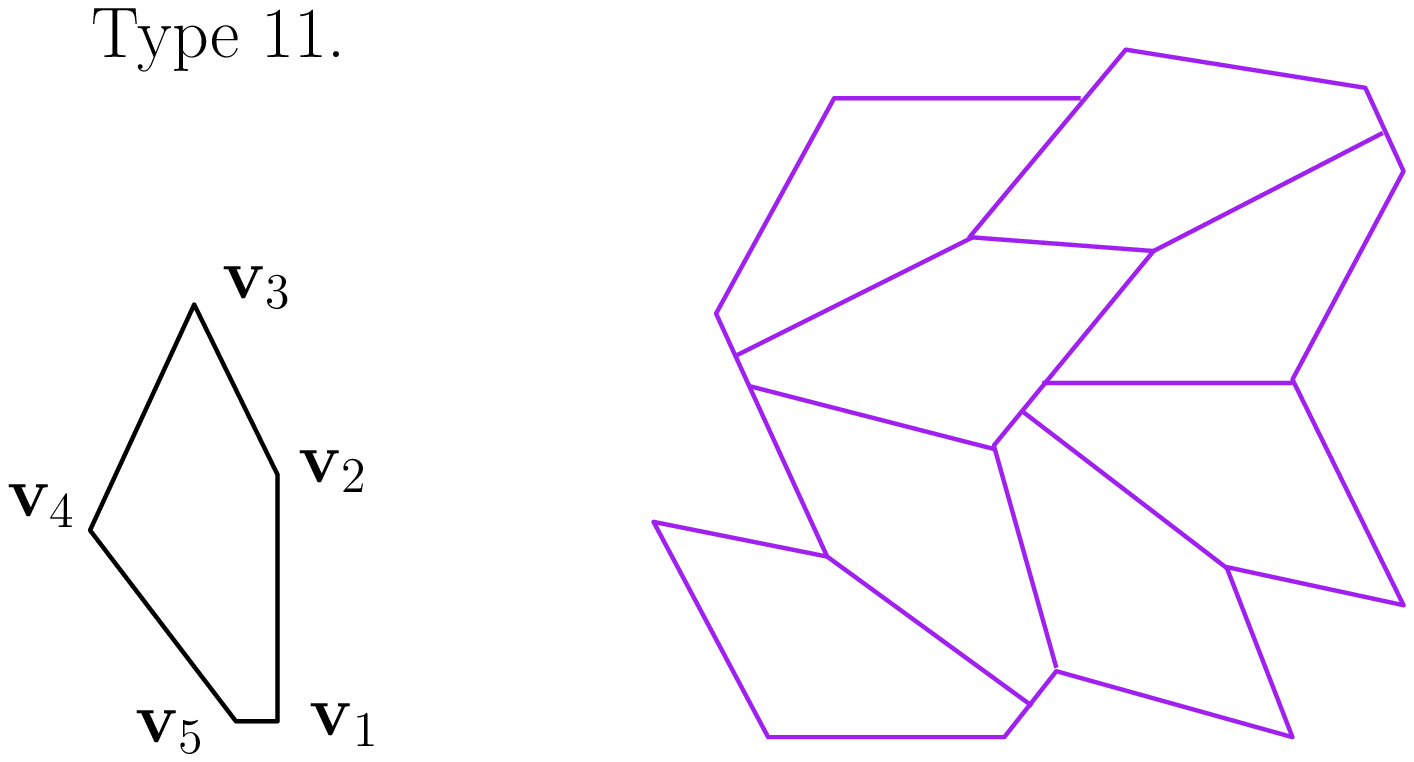}

\vspace{0.5cm}
\centering
\includegraphics[scale=0.45]{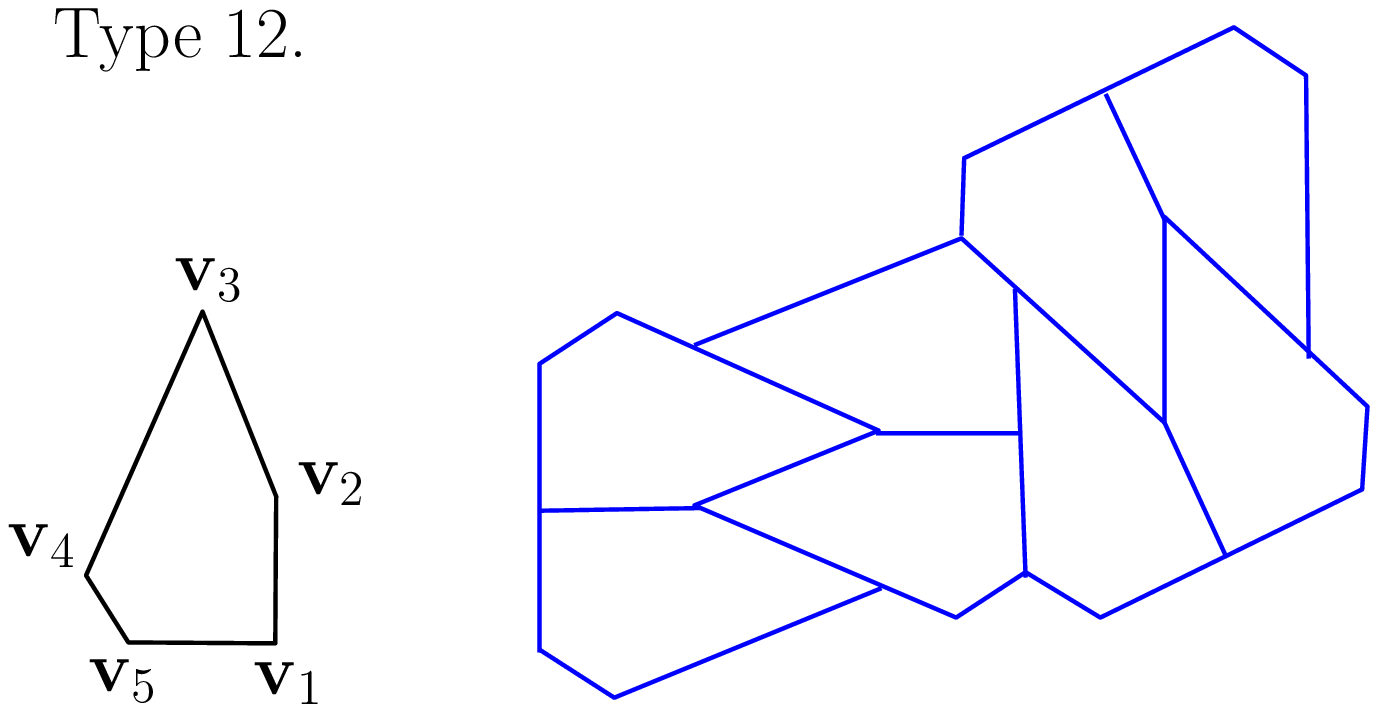}

\vspace{0.45cm}
\centering
\includegraphics[scale=0.45]{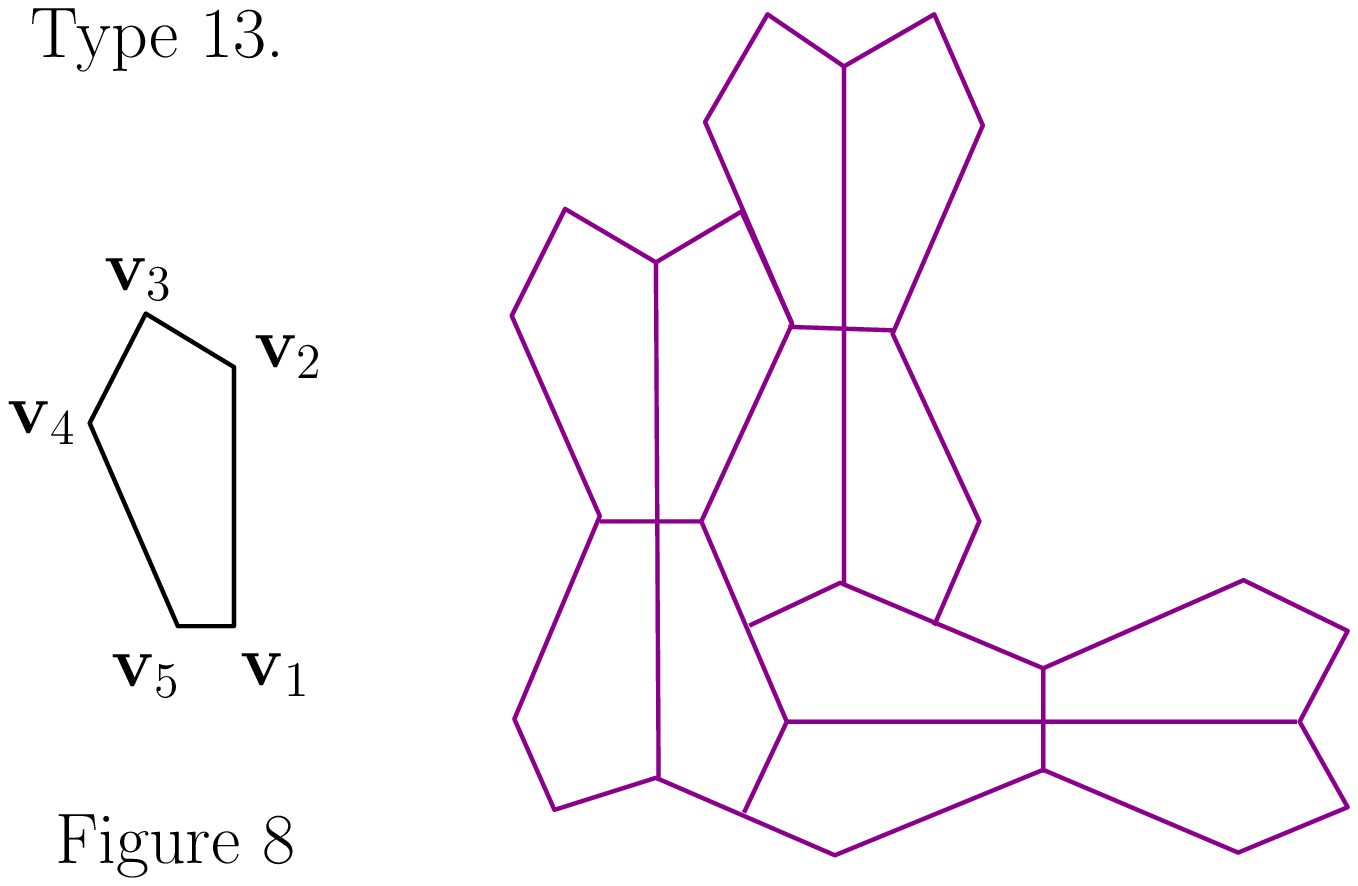}
\end{figure}

Marjorie Rice died on July 2, 2017, at the age of 94. A lobby floor of the Mathematical Association of America in Washington is paved by one of Rice's pentagon tile in her honor. On July 11, 2017, Quanta Magazine \cite{wolc} published a article in her memory.

Rice's method was systematic, in the sense based on some geometric principle. Nevertheless, it was not able to guarantee the completeness of the list. In 1985, Rolf Stein reported another one.

\medskip\noindent
{\bf Theorem 7 (Stein \cite{stei85}).} {\it A convex pentagon $P_5$ can tile the whole plane $\mathbb{E}^2$ if it satisfies the following conditions:

\noindent\begin{enumerate}
\item[\bf (14).] $\alpha_1={\pi\over 2}$, $2\alpha_2+\alpha_3=2\pi $, $\alpha_3+\alpha_5=\pi$, and $2\ell_1=2\ell_3=\ell_4=\ell_5$.
\end{enumerate}}

\medskip
\begin{figure}[!ht]
\centering
\includegraphics[scale=0.45]{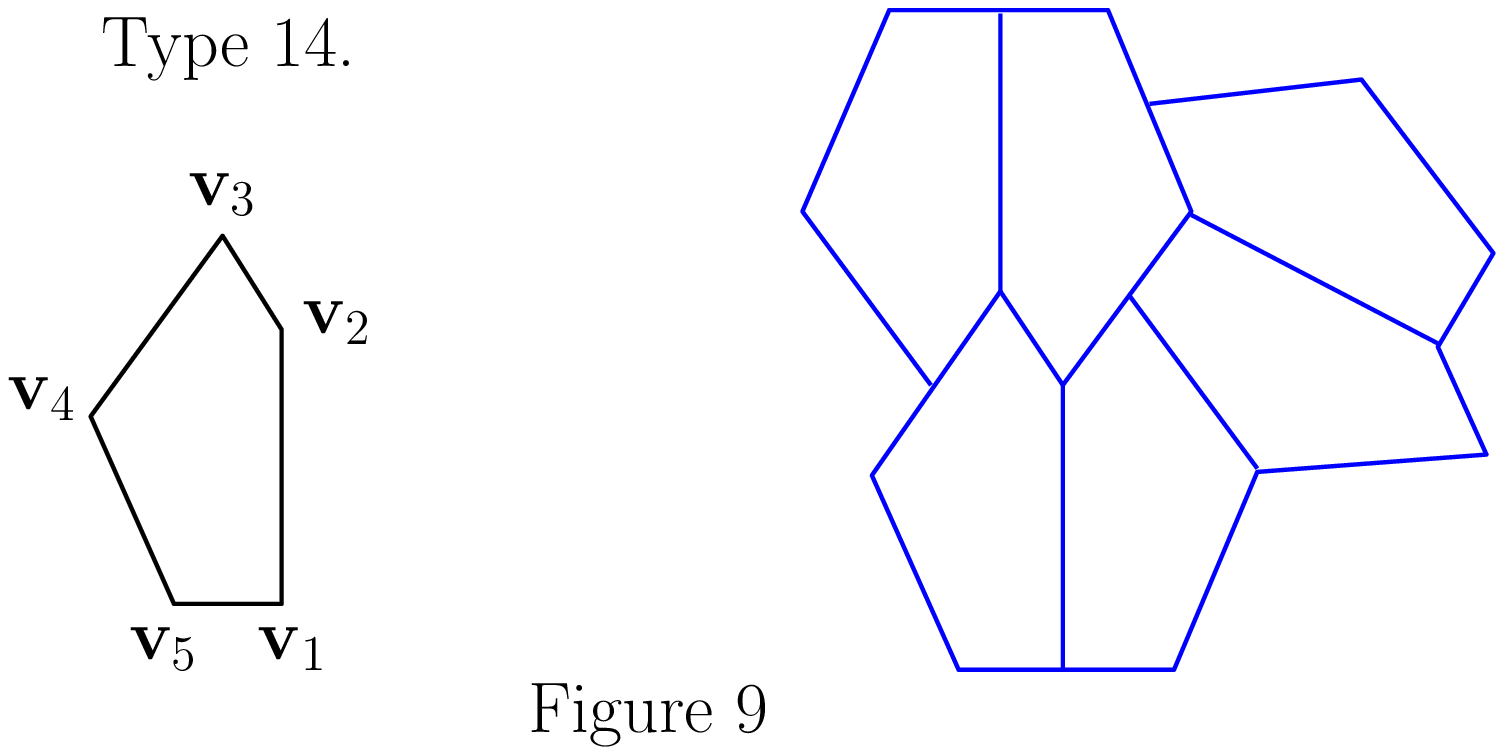}
\end{figure}

\smallskip
Stein's discovery is illustrated by Figure 9.

\vspace{0.6cm}
\noindent
{\Large\bf 4. Fifteen, And Only Fifteen}

\bigskip\noindent
Let $\mathcal{T}$ denote a tiling of $\mathbb{E}^2$ with congruent tiles. A {\it symmetry} of $\mathcal{T}$ is an isometry of $\mathbb{E}^2$ that maps the tiles of $\mathcal{T}$ onto tiles of $\mathcal{T}$, and the {\it symmetry group} $\mathcal{G}$ of $\mathcal{T}$ is the collection of all such symmetries associated with isometry multiplications. Two tiles $T_1$ and $T_2$ of $\mathcal{T}$ are said to be equivalent if there is a symmetry $\sigma \in \mathcal{G}$ such that $\sigma (T_1)=T_2$. If all the tiles of $\mathcal{T}$ are equivalent to one tile $T$, the tiling $\mathcal{T}$ is said to be {\it transitive} and $T$ is called a {\it transitive tile}. Then, Hilbert's problem can be reformulated as:

\medskip
{\it Is every polytope which can tile the whole space a transitive tile}?

\medskip
A tiling $\mathcal{T}$ of $\mathbb{E}^2$ by identical convex pentagons is called an {\it $n$-block transitive tiling} if it has a block $B$ consists of $n$ (minimum) connected tiles such that $\mathcal{T}$ is a transitive tiling of $B$. If a convex pentagon $T$ can form a tiling is as small as $n$-block transitive, we call it an {\it $n$-block transitive tile}. Clearly, all the tiles of Types 1-5 are one-block transitive. In other words, they are transitive tiles. According to \cite{mann18,scha78}, all the tiles of Types 5-14 except Type 9 are two-block transitive and the tiles of Type 9 are three-block transitive.

From the intuitive point of view, it is reasonable to believe that periodic structure is inevitable in pentagon tilings and the period can not be too large. Based on this belief, Mann, McLoud-Mann and Von Derau \cite{mann18} developed an algorithm for enumerating all the $n$-block transitive pentagon tiles. When they check the three-block case, surprisingly, a new type of pentagon tiles is discovered.

\medskip\noindent
{\bf Theorem 8 (Mann, McLoud-Mann and Von Derau \cite{mann18}).} {\it A convex pentagon $P_5$ can tile the whole plane $\mathbb{E}^2$ if it satisfies the following conditions:

\noindent\begin{enumerate}
\item[\bf (15).] $\alpha_1={\pi\over 3}$, $\alpha_2={{3\pi}\over 3}$, $\alpha_3={{7\pi}\over {12}}$, $\alpha_4={{\pi}\over {2}}$, $\alpha_5={{5\pi}\over {6}}$, and $\ell_1=2\ell_2=2\ell_4=2\ell_5$.
\end{enumerate}}

\medskip
\begin{figure}[!ht]
\centering
\includegraphics[scale=0.39]{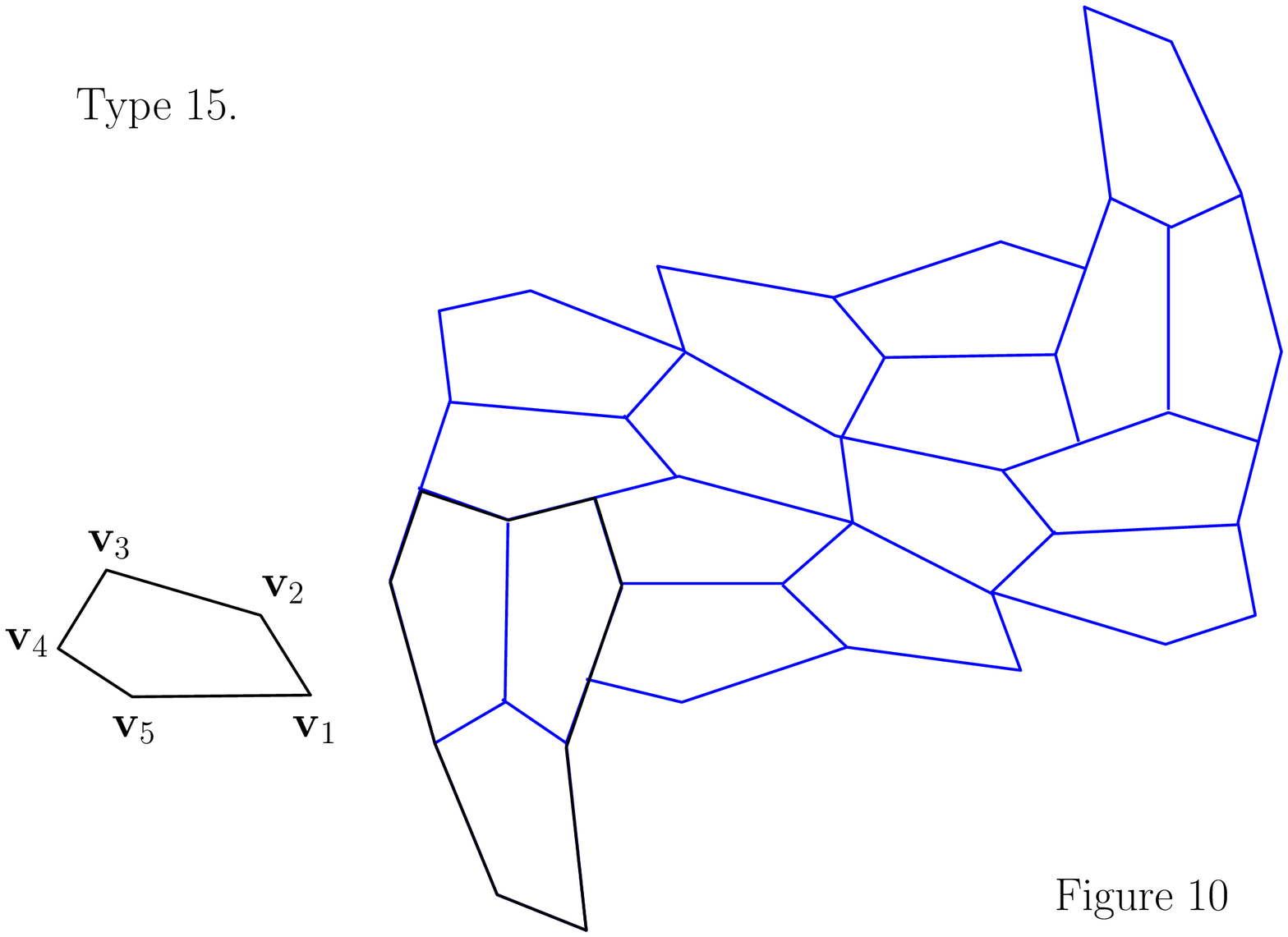}
\end{figure}

\medskip
Figure 10 shows a tiling pattern of their discovery.

\medskip\noindent
{\bf Remark 2.} It was shown by Mann, McLoud-Mann and Von Derau \cite{mann18} that there is no other $n$-block transitive pentagon tile with $n\le 4$. Then, the completeness of the list emerges again.

\medskip
Since Hales' computer proof for the Kepler conjecture, more and more geometers turn to computer for help when their mathematical problems can be reduced into a big number of cases. Characterizing all the pentagon tiles seems to be a perfect candidate for such purpose. Of course, one should first get round the periodic assumption for pentagon tilings.

In 2017, one century after Bieberbach proposed the characterization problem, Micha\"el Rao announced a computer proof for the completeness of the known pentagon tile list. Rao's approach is based on graph expression. First he proved that, if a pentagon tiles the plane, then it can form a tiling
such that every vertex type has positive density. Clearly, this is a weak version of the periodic tiling. Second, it was shown that there are only finite number of possible vertex types in the modified pentagon tiling. In fact, he reduced them into 371 types. Then, by testing the 371 cases Rao proved the following theorem.

\medskip\noindent
{\bf Theorem 9 (Rao \cite{rao17}).} {\it A convex pentagon $P_5$ can tile the whole plane if and only if it belongs to one of the fifteen types listed in Theorems {\rm 3--8}.}

\medskip
Computer proofs are still not as acceptable as transparent logical proofs within the mathematical community. However, we have to admit that the complexities of mathematical problems changes from zero to infinity, and perhaps there indeed exist problems which have no transparent logical proofs.

\vspace{0.6cm}
\noindent
{\Large\bf 5. Multiple Tilings, Basic}

\bigskip\noindent
Let $K$ denote an $n$-dimensional convex body with interior ${\rm int}(K)$, boundary $\partial (K)$, and volume ${\rm vol}(K)$. In particular, let $D$ denote a two-dimensional convex domain.

Assume that $\mathcal{F}=\{ K_1, K_2, K_3, \ldots \}$ is a family of convex bodies in $\mathbb{E}^n$ and $k$ is a positive integer. We call $\mathcal{F}$ a {\it $k$-fold tiling} of $\mathbb{E}^n$ if every point ${\bf x}\in \mathbb{E}^n$ belongs to at least $k$ of these convex bodies
and every point ${\bf x}\in \mathbb{E}^n$ belongs to at most $k$ of the ${\rm int}(K_i)$. In other words, a $k$-fold tiling of $\mathbb{E}^n$ is both a {\it $k$-fold packing} and a {\it $k$-fold covering} in $\mathbb{E}^n$. In particular, we call a $k$-fold tiling of $\mathbb{E}^n$ a {\it $k$-fold congruent tiling}, a {\it $k$-fold translative tiling}, or a {\it $k$-fold lattice tiling} if all $K_i$ are congruent to $K_1$, all $K_i$ are translates of $K_1$, or all $K_i$ are translates of $K_1$ and the translative vectors form a lattice in $\mathbb{E}^n$, respectively. In these particular cases, we call $K_1$ a {\it $k$-fold congruent tile}, a {\it $k$-fold translative tile} or a {\it $k$-fold lattice tile}, respectively.

For a fixed convex body $K$, we define $\tau^\bullet (K)$ to be the smallest integer $k$ such that $K$ can form a $k$-fold congruent tiling in $\mathbb{E}^n$, $\tau (K)$ to be the smallest integer $k$ such that $K$ can form a $k$-fold translative tiling in $\mathbb{E}^n$,
and $\tau^*(K)$ to be the smallest integer $k$ such that $K$ can form a $k$-fold lattice tiling in $\mathbb{E}^n$. For convenience, if $K$ can not form any multiple congruent tiling, translative tiling or lattice tiling, we will define $\tau^\bullet(K)=\infty$, $\tau (K)=\infty$ or $\tau^*(K)=\infty$, respectively. Clearly, for every convex body $K$ we have
$$\tau^\bullet(K)\le \tau (K)\le \tau^*(K).$$

By looking at the separating hyperplanes between tangent neighbours, it is obvious that a convex body can form a multiple tiling only if it is a polytope.

If $\sigma $ is a non-singular affine linear transformation from $\mathbb{E}^n$ to $\mathbb{E}^n$, then $\mathcal{F}=\{ K_1, K_2, K_3, \ldots \}$ forms a $k$-fold tiling of $\mathbb{E}^n$ if and only if $\mathcal{F}'=\{ \sigma (K_1), \sigma (K_2),$ $\sigma (K_3),$ $\ldots \}$ forms a $k$-fold tiling of $\mathbb{E}^n$. Consequently, for any $n$-dimensional convex body $K$ and any non-singular affine linear transformation $\sigma$ we have
both
$$\tau (\sigma (K))=\tau (K)$$
and
$$\tau^* (\sigma (K))=\tau^* (K).$$
Unfortunately, $\tau^\bullet (K)$ is not an invariant for the linear transformation group.

Clearly, one-fold tilings are the usual tilings. In the plane, we have
$$\tau (D)=\tau^*(D)=1$$ if and only if $D$ is a parallelogram or a centrally symmetric hexagon, and
$$\tau^\bullet (D)=1$$ if and only if $D$ is a triangle, a quadrilateral, a pentagon belonging to one of the fifteen types listed in Theorem 3--8,
or an hexagon belonging to one of the three types listed in Theorem 2.

Since 1936, multiple tilings has been studied by Furtw\"angler \cite{furt}, Haj\'os \cite{hajo}, Robinson \cite{robi}, Bolle \cite{boll94}, Gravin, Robins and Shiryaev \cite{grs} and many others. Nevertheless, many natural problems are still open. In the forthcoming sections we will introduce some fascinating new results about multiple tilings in the plane.

\vspace{0.6cm}
\noindent
{\Large\bf 6. Multiple Lattice Tilings}

\bigskip\noindent
In 1994, Bolle studied the two-dimensional lattice multiple tilings. He proved the following criterion:

\medskip\noindent
{\bf Lemma 2 (Bolle \cite{boll94}).} {\it A convex polygon is a $k$-fold lattice tile for a lattice $\Lambda$ and some positive integer $k$ if and only if the following conditions are satisfied:

\noindent
{\bf 1.} It is centrally symmetric.

\noindent
{\bf 2.} When it is centered at the origin, in the relative interior of each edge $G$ there is a point of ${1\over 2}\Lambda $.

\noindent
{\bf 3.} If the midpoint of $G$ is not in ${1\over 2}\Lambda $, then $G$ is a lattice vector of $\Lambda $.}

\medskip
Based on Bolle's criterion, Gravin, Robins and Shiryaev \cite{grs} discovered the following example.

\medskip\noindent
{\bf Example 1.} Let $\Lambda $ denote the two-dimensional integer lattice $\mathbb{Z}^2$ and let $D_8$ denote the polygon with vertices
${\bf v}_1=({1\over 2},-{3\over 2})$, ${\bf v}_2=({3\over 2},-{1\over 2})$, ${\bf v}_3=({3\over 2},{1\over 2})$, ${\bf v}_4=({1\over 2},{3\over 2})$, ${\bf v}_5=-{\bf v}_1$, ${\bf v}_6=-{\bf v}_2$, ${\bf v}_7=-{\bf v}_3$ and ${\bf v}_8=-{\bf v}_4$, as shown by Figure 11. Then $D_8+\Lambda $ is a seven-fold lattice tiling of $\mathbb{E}^2$.

\begin{figure}[!ht]
\centering
\includegraphics[scale=0.5]{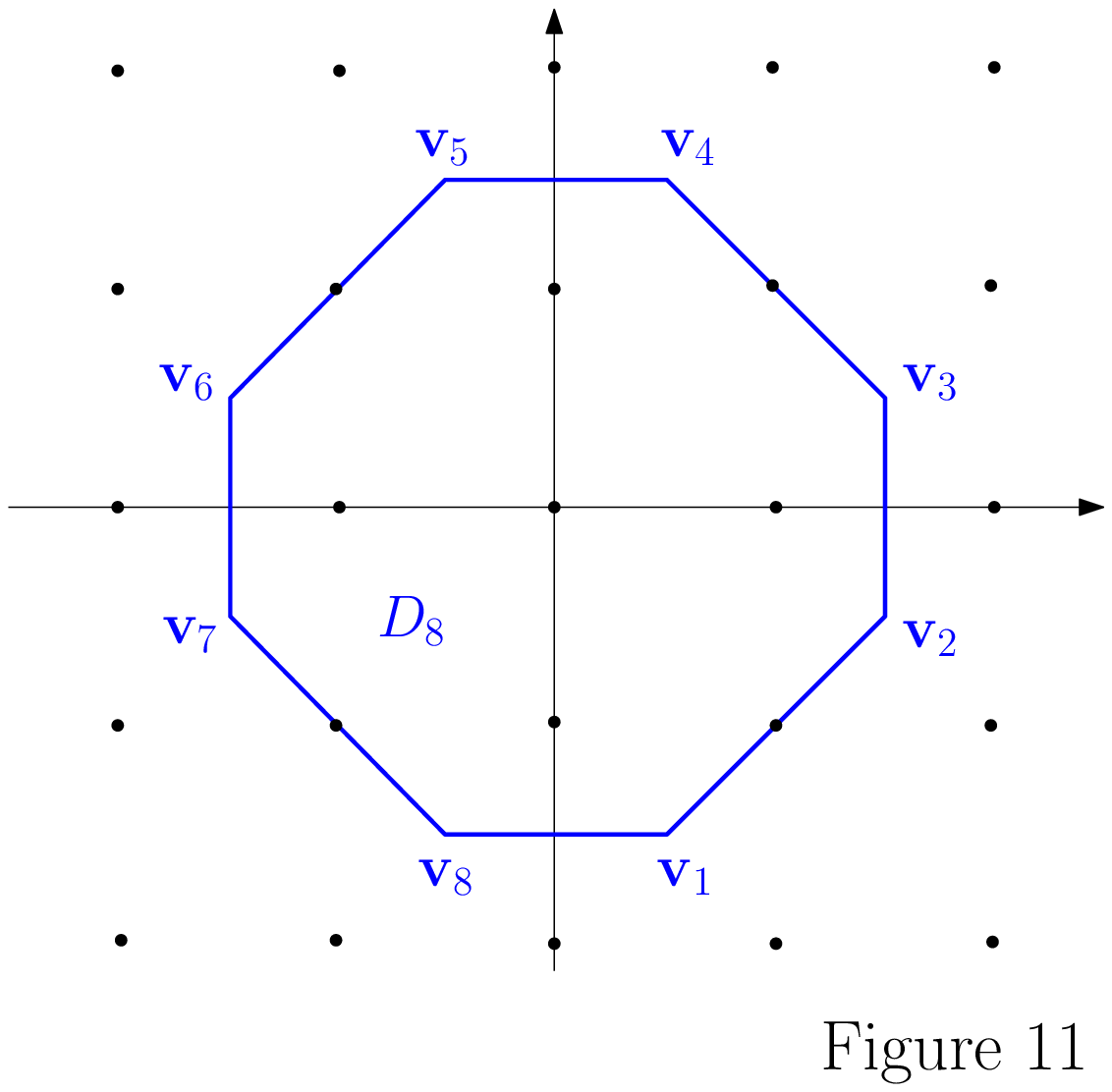}
\end{figure}

Let $\mathcal{D}$ denote the family of all two-dimensional convex domains and let $\mathcal{P}_{2m}$ denote the family of all centrally symmetric
$2m$-gons. Example 1 is a strong evidence for conjecturing
$$\min_{D\in \mathcal{D}\setminus \{\mathcal{P}_4\cup \mathcal{P}_6\}} \tau^*(D)\ge 7.\eqno(1)$$

A multiple tiling is both a multiple packing and a multiple covering. In the literature, multiple packing has been studied by many authors (see Zong \cite{zong14}). Let $\delta_k(D)$ denote the density of the densest $k$-fold lattice packing of $D$ in $\mathbb{E}^2$. In particular, we have the following lemma.

\medskip\noindent
{\bf Lemma 3 (Dumir and Hans-Gill \cite{dumi}, G. Fejes T\'oth \cite{feje}).} {\it If $k=2,$ $3$ or $4$, then
$$\delta_k(D)=k\cdot \delta_1 (D)$$
holds for every two-dimensional centrally symmetric convex domain $D$.}

\medskip
In 2017, Yang and Zong obtained the following unexpected result.

\medskip\noindent
{\bf Theorem 10 (Yang and Zong \cite{yz1}).} {\it If $D$ is a two-dimensional convex domain which is neither a parallelogram nor a centrally symmetric hexagon, then we have
$$\tau ^*(D)\ge 5,$$
where the equality holds at some particular decagons.}

\medskip
This theorem can be easily deduced from Lemma 3 and Lemma 2. One of the interesting things is to notice that (1) is not true. Let $\Lambda $ to be the integer lattice $\mathbb{Z}^2$ and let $D_{10}$ denote the decagon with ${\bf u}_1=(0, 1)$, ${\bf u}_2=(1, 1)$, ${\bf u}_3=({3\over 2}, {1\over 2})$, ${\bf u}_4=({3\over 2}, 0)$, ${\bf u}_5=( 1,-{1\over 2})$, ${\bf u}_6=-{\bf u}_1$, ${\bf u}_7=-{\bf u}_2$, ${\bf u}_8=-{\bf u}_3$, ${\bf u}_9=-{\bf u}_4$ and ${\bf u}_{10}=-{\bf u}_5$ as the middle points of its edges, as shown by Figure 12. By Lemma 2, it can be easily verified that $D_{10}+\Lambda $ is indeed a five-fold lattice tiling of $\mathbb{E}^2$.

\medskip
\begin{figure}[!ht]
\centering
\includegraphics[scale=0.44]{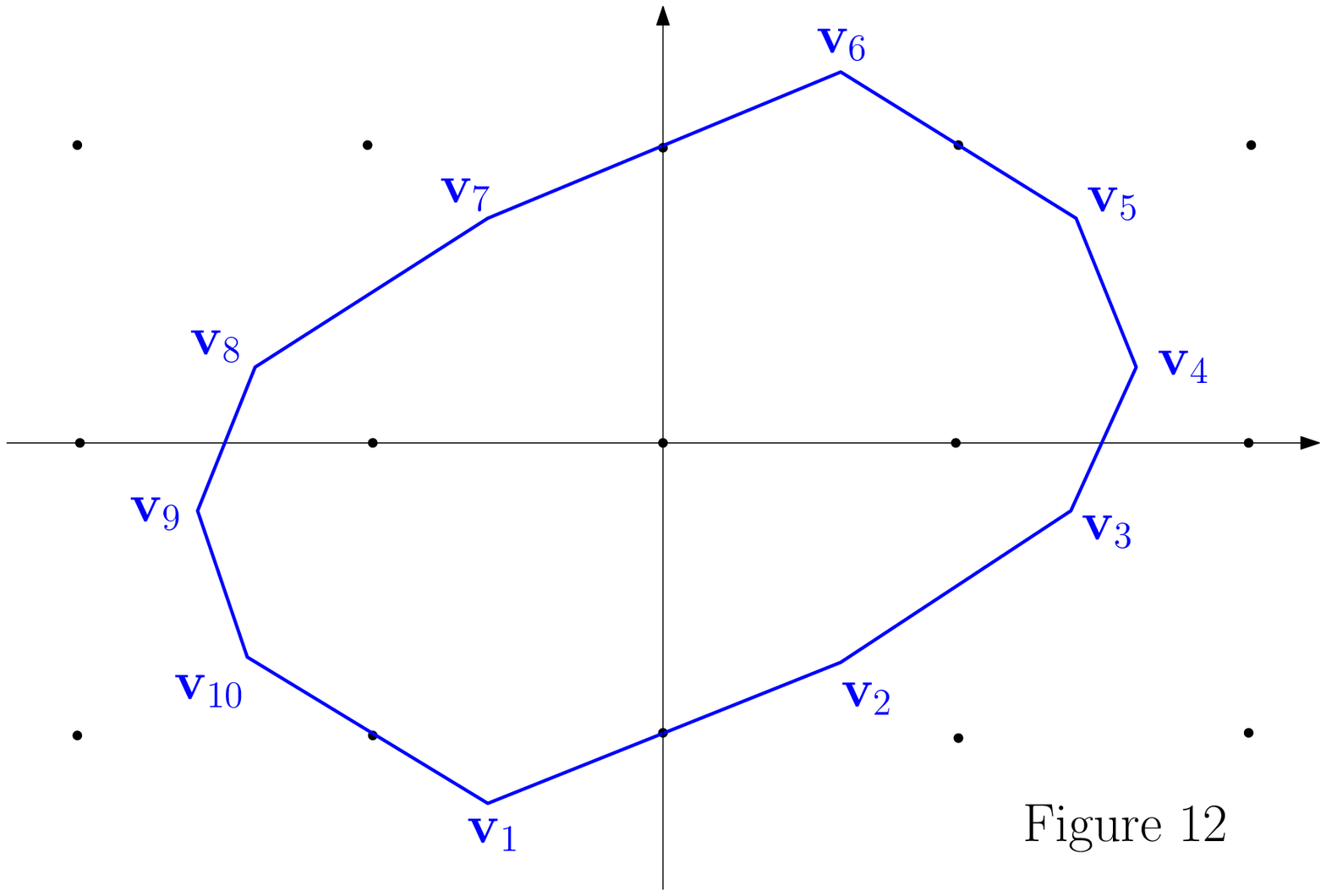}
\end{figure}

Even more unexpected, by studying lattice polygons, all the five-fold lattice tiles can be nicely characterized. They are two classes octagons and one class of decagons, besides the parallelograms and the centrally symmetric hexagons.

\medskip\noindent
{\bf Theorem 11 (Zong \cite{zong}).} {\it A convex domain can form a five-fold lattice tiling of the Euclidean plane if and only if it is a parallelogram, a centrally symmetric hexagon, under a suitable affine linear transformation, a centrally symmetric octagon with vertices
${\bf v}_1=(-\alpha , -{3\over 2})$, ${\bf v}_2=(1-\alpha , -{3\over 2})$, ${\bf v}_3=(1+\alpha , -{1\over 2})$, ${\bf v}_4=(1-\alpha , {1\over 2})$, ${\bf v}_5=-{\bf v}_1$, ${\bf v}_6=-{\bf v}_2$, ${\bf v}_7=-{\bf v}_3$ and ${\bf v}_8=-{\bf v}_4$, where $0<\alpha <{1\over 4},$ or with vertices
${\bf v}_1=(\beta , -2)$, ${\bf v}_2=(1+\beta , -2)$, ${\bf v}_3=(1-\beta , 0)$, ${\bf v}_4=(\beta , 1)$, ${\bf v}_5=-{\bf v}_1$, ${\bf v}_6=-{\bf v}_2$, ${\bf v}_7=-{\bf v}_3$, ${\bf v}_8=-{\bf v}_4$, where ${1\over 4}<\beta <{1\over 3}$, or a centrally symmetric decagon with ${\bf u}_1=(0, 1)$, ${\bf u}_2=(1, 1)$, ${\bf u}_3=({3\over 2}, {1\over 2})$, ${\bf u}_4=({3\over 2}, 0)$, ${\bf u}_5=( 1,-{1\over 2})$, ${\bf u}_6=-{\bf u}_1$, ${\bf u}_7=-{\bf u}_2$, ${\bf u}_8=-{\bf u}_3$, ${\bf u}_9=-{\bf u}_4$ and ${\bf u}_{10}=-{\bf u}_5$ as the middle points of its edges.}

\medskip
Needless to say that the proof for this theorem is very complicated and technical. In fact, in the course to discover and to prove this theorem, several mistakes were made. Fortunately, all the detected errors could be corrected.

Let $D_8(\alpha )$ denote the octagon with vertices ${\bf v}_1=(-\alpha , -{3\over 2})$, ${\bf v}_2=(1-\alpha , -{3\over 2})$, ${\bf v}_3=(1+\alpha , -{1\over 2})$, ${\bf v}_4=(1-\alpha , {1\over 2})$, ${\bf v}_5=-{\bf v}_1$, ${\bf v}_6=-{\bf v}_2$, ${\bf v}_7=-{\bf v}_3$ and ${\bf v}_8=-{\bf v}_4$, where $0<\alpha <{1\over 4},$ and let $\Lambda $ denote the integer lattice, as shown by Figure 13. By Lemma 2 it can be verified that $D_8(\alpha )+\Lambda $ is indeed a five-fold lattice tiling.

\begin{figure}[!ht]
\centering
\includegraphics[scale=0.5]{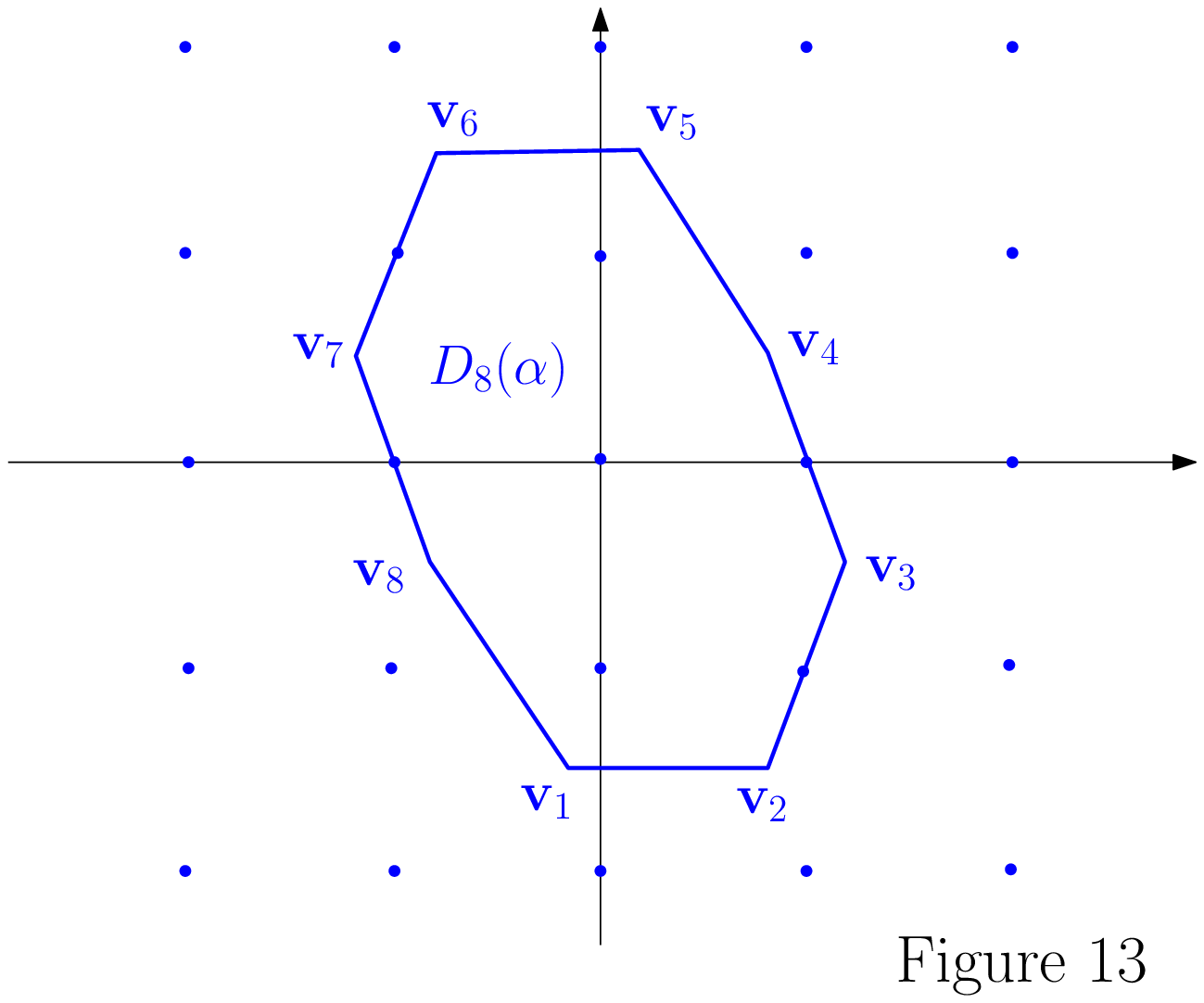}
\end{figure}

Similarly, let $D_8(\beta )$ denote the octagon with vertices ${\bf v}_1=(\beta , -2)$, ${\bf v}_2=(1+\beta , -2)$, ${\bf v}_3=(1-\beta , 0)$, ${\bf v}_4=(\beta , 1)$, ${\bf v}_5=-{\bf v}_1$, ${\bf v}_6=-{\bf v}_2$, ${\bf v}_7=-{\bf v}_3$, ${\bf v}_8=-{\bf v}_4$, where ${1\over 4}<\beta <{1\over 3}$, and let $\Lambda $ denote the integer lattice, as shown by Figure 14. By Lemma 2 it can be verified that $D_8(\beta )+\Lambda $ is indeed a five-fold lattice tiling as well.

\begin{figure}[!ht]
\includegraphics[scale=0.5]{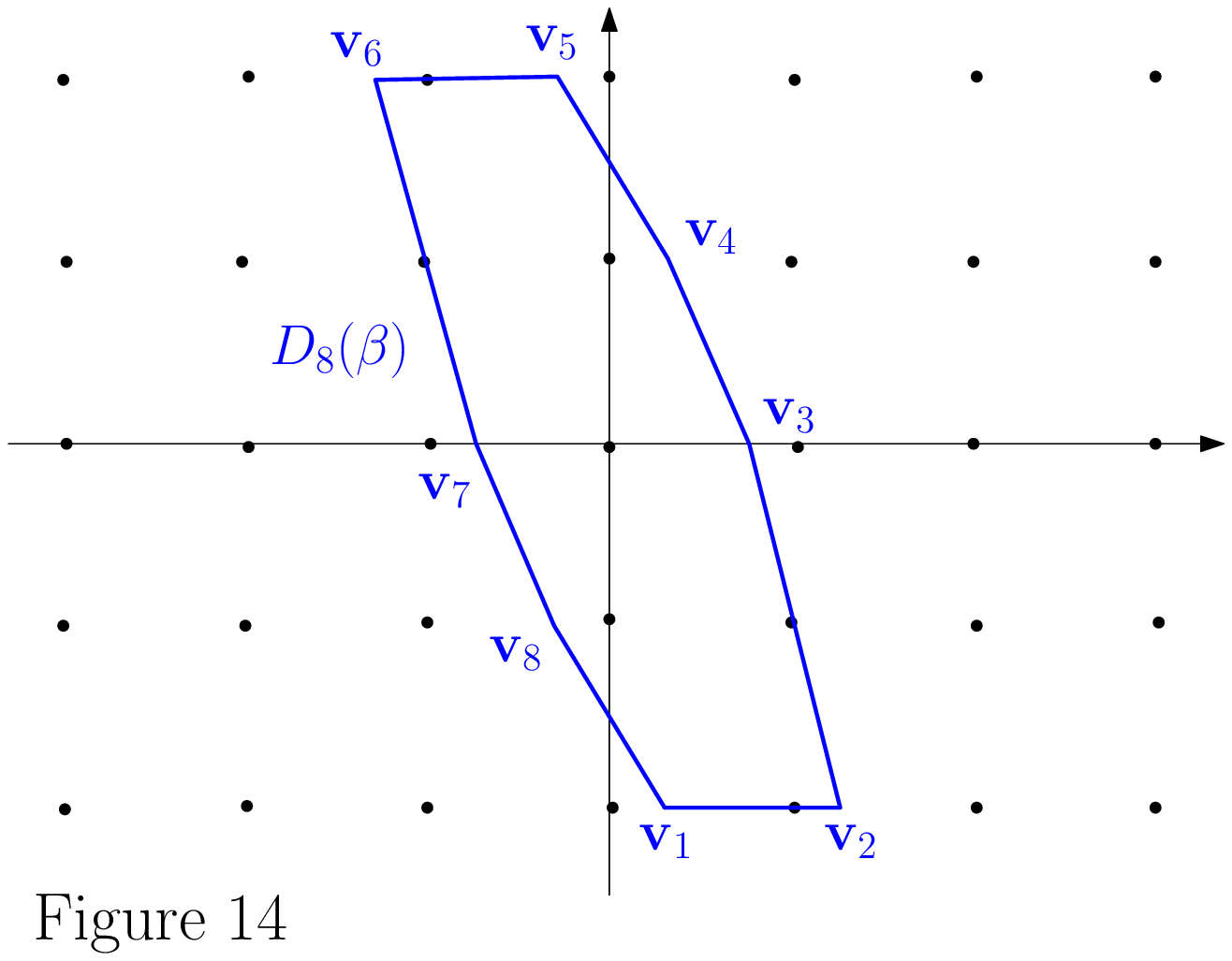}
\end{figure}

\medskip\noindent
{\bf Theorem 12 (Zong \cite{zong}).} {\it Let $W$ denote the quadrilateral with vertices ${\bf w}_1=(-{1\over 2}, 1)$, ${\bf w}_2=(-{1\over 2}, {3\over 4})$, ${\bf w}_3=(-{2\over 3}, {2\over 3})$ and ${\bf w}_4=(-{3\over 4}, {3\over 4})$. A centrally symmetric convex decagon with ${\bf u}_1=(0,1)$, ${\bf u}_2=(1,1)$, ${\bf u}_3=({3\over 2}, {1\over 2})$, ${\bf u}_4=({3\over 2}, 0)$, ${\bf u}_5=(1,-{1\over 2})$, ${\bf u}_6=-{\bf u}_1$, ${\bf u}_7=-{\bf u}_2$, ${\bf u}_8=-{\bf u}_3$, ${\bf u}_9=-{\bf u}_4$ and ${\bf u}_{10}=-{\bf u}_5$ as the middle points of its edges if and only if one of its vertices is an interior point of $W$.}

\vspace{0.6cm}
\noindent
{\Large\bf 7. Multiple Translative Tilings}

\bigskip\noindent
In 2012, Gravin, Robins and Shiryaev \cite{grs} proved that an $n$-dimensional convex body can form a multiple translative tiling of the space only if it is a centrally symmetric polytope with centrally symmetric facets. Therefore, in the plane we only need to deal with the centrally symmetric polygons.

Let $P_{2m}$ denote a centrally symmetric convex $2m$-gon centered at the origin, with vertices ${\bf v}_1$, ${\bf v}_2$, $\ldots $, ${\bf v}_{2m}$
enumerated in the clock-order. Let $G_i$ denote the edge with ends ${\bf v}_i$ and ${\bf v}_{i+1}$, and write $V=\{ {\bf v}_1, {\bf v}_2, \ldots , {\bf v}_{2m}\}$.

Assume that $P_{2m}+X$ is a $\tau (P_{2m})$-fold translative tiling in $\mathbb{E}^2$, where $X=\{{\bf x}_1, {\bf x}_2, {\bf x}_3, \ldots \}$ is a discrete multiset with ${\bf x}_1={\bf o}$. By observing the local structure of $P_{2m}+X$ at the vertices ${\bf v}\in V+X$, Yang and Zong \cite{yz2,yz3} discovered some fascinating results.

Let $X^{\bf v}$ denote the subset of $X$ consisting of all points ${\bf x}_i$ such that
$${\bf v}\in \partial (P_{2m})+{\bf x}_i.$$
Since $P_{2m}+X$ is a multiple tiling, the set $X^{\bf v}$ can be divided into disjoint subsets $X^{\bf v}_1$, $X^{\bf v}_2$, $\ldots ,$ $X^{\bf v}_t$ such that the translates in $P_{2m}+X^{\bf v}_j$ can be re-enumerated as $P_{2m}+{\bf x}^j_1$, $P_{2m}+{\bf x}^j_2$, $\ldots $, $P_{2m}+{\bf x}^j_{s_j}$ satisfying the following conditions:

\smallskip
\noindent
{\bf 1.} {\it ${\bf v}\in \partial (P_{2m})+{\bf x}^j_i$ holds for all $i=1, 2, \ldots, s_j.$}

\noindent
{\bf 2.} {\it Let $\angle^j_i$ denote the inner angle of $P_{2m}+{\bf x}^j_i$ at ${\bf v}$ with two half-line edges $L^j_{i,1}$ and $L^j_{i,2}$ such that $L^j_{i,1}$, ${\bf x}^j_i-{\bf v}$ and $L^j_{i,2}$ are in clock order. Then, the inner angles join properly as
$$L^j_{i,2}=L^j_{i+1,1}$$
holds for all $i=1,$ $2,$ $\ldots ,$ $s_j$, where $L^j_{s_j+1,1}=L^j_{1,1}$.}

\smallskip
For convenience, we call such a sequence $P_{2m}+{\bf x}^j_1$, $P_{2m}+{\bf x}^j_2$, $\ldots $, $P_{2m}+{\bf x}^j_{s_j}$ an {\it adjacent wheel} at ${\bf v}$. It is easy to see that
$$\sum_{i=1}^{s_j}\angle^j_i =2w_j\cdot \pi$$
hold for positive integers $w_j$. Then we define
$$\phi ({\bf v})=\sum_{j=1}^tw_j= {1\over {2\pi }}\sum_{j=1}^t\sum_{i=1}^{s_j}\angle^j_i$$
and
$$\varphi ({\bf v})=\sharp \left\{ {\bf x}_i:\ {\bf x}_i\in X,\ {\bf v}\in {\rm int}(P_{2m})+{\bf x}_i\right\}.$$

Clearly, if $P_{2m}+X$ is a $\tau (P_{2m})$-fold translative tiling of $\mathbb{E}^2$, then
$$\tau (P_{2m})= \varphi ({\bf v})+\phi ({\bf v})\eqno(2)$$
holds for all ${\bf v}\in V+X$. By detailed analysis based on (2), Yang and Zong \cite{yz2,yz3} obtained the following results.

\medskip\noindent
{\bf Theorem 13 (Yang and Zong \cite{yz2}).} {\it If $D$ is a two-dimensional convex domain which is neither a parallelogram nor a centrally symmetric hexagon, then we have
$$\tau (D)\ge 5,$$
where the equality holds if $D$ is some particular centrally symmetric octagon or some particular centrally symmetric decagon.}

\medskip\noindent
{\bf Remark 3.} It is known that
$$\tau (D)\le \tau^*(D)$$
holds for every convex domain $D$. Therefore, Theorem 13 implies Theorem 10.

\smallskip
At this point, it should be interesting to characterize all the five-fold translative tiles, in particular to figure out if they are the known
five-fold lattice tiles.

\medskip\noindent
{\bf Theorem 14 (Yang and Zong \cite{yz3}).} {\it A convex domain can form a five-fold translative tiling of the Euclidean plane if and only if it is a parallelogram, a centrally symmetric hexagon, under a suitable affine linear transformation, a centrally symmetric octagon with vertices
${\bf v}_1=\left( {3\over 2}-{{5\alpha }\over 4}, -2\right)$, ${\bf v}_2=\left( -{1\over 2}-{{5\alpha }\over 4}, -2\right)$, ${\bf v}_3=\left( {{\alpha }\over 4}-{3\over 2}, 0\right)$, ${\bf v}_4=\left( {{\alpha }\over 4}-{3\over 2}, 1\right)$, ${\bf v}_5=-{\bf v}_1$, ${\bf v}_6=-{\bf v}_2$, ${\bf v}_7=-{\bf v}_3$ and ${\bf v}_8=-{\bf v}_4$, where $0<\alpha <{2\over 3}$, or with vertices ${\bf v}_1=(2-\beta, -3),$ ${\bf v}_2=(-\beta, -3),$ ${\bf v}_3=(-2, -1),$ ${\bf v}_4=(-2, 1),$ ${\bf v}_5=-{\bf v}_1$, ${\bf v}_6=-{\bf v}_2$, ${\bf v}_7=-{\bf v}_3$ and ${\bf v}_8=-{\bf v}_4$, where $0<\beta \le 1$, or a centrally symmetric decagon with $\mathbf{u}_{1}=(0,1),$ $\mathbf{u}_{2}=(1,1),$ $\mathbf{u}_{3}=({3\over 2}, {1\over 2}),$ $\mathbf{u}_{4}=({3\over 2},0),$ $\mathbf{u}_{5}=(1,-{1\over 2}),$ $\mathbf{u}_{6}=-\mathbf{u}_{1},$ $\mathbf{u}_{7}=-\mathbf{u}_{2},$ $\mathbf{u}_{8}=-\mathbf{u}_{3},$ $\mathbf{u}_{9}=-\mathbf{u}_{4}$ and $\mathbf{u}_{10}=-\mathbf{u}_{5}$ as the middle points of its edges.}

\medskip
Needless to say, the proofs for Theorem 13 and Theorem 14 are extremely complicated. However, the shapes of the resulted polygons are quite elegant. The decagons in Theorem 14 are the decagons in Theorem 11, which were shown by Figure 12.

Let $D'_8(\alpha )$ denote the octagon with vertices ${\bf v}_1=\left( {3\over 2}-{{5\alpha }\over 4}, -2\right)$, ${\bf v}_2=\left( -{1\over 2}-{{5\alpha }\over 4}, -2\right)$, ${\bf v}_3=\left( {{\alpha }\over 4}-{3\over 2}, 0\right)$, ${\bf v}_4=\left( {{\alpha }\over 4}-{3\over 2}, 1\right)$, ${\bf v}_5=-{\bf v}_1$, ${\bf v}_6=-{\bf v}_2$, ${\bf v}_7=-{\bf v}_3$ and ${\bf v}_8=-{\bf v}_4$, where $0<\alpha <{2\over 3}$, and let
$\Lambda (\alpha)$ denote the lattice generated by ${\bf u}_1=(2,0)$ and ${\bf u}_2=(1+{\alpha \over 2}, 1)$, as shown by Figure 15. By Lemma 2 it can be verified that $D'_8(\alpha )+\Lambda (\alpha )$ is indeed a five-fold translative tiling in $\mathbb{E}^2$.

\begin{figure}[!ht]
\centering
\includegraphics[scale=0.55]{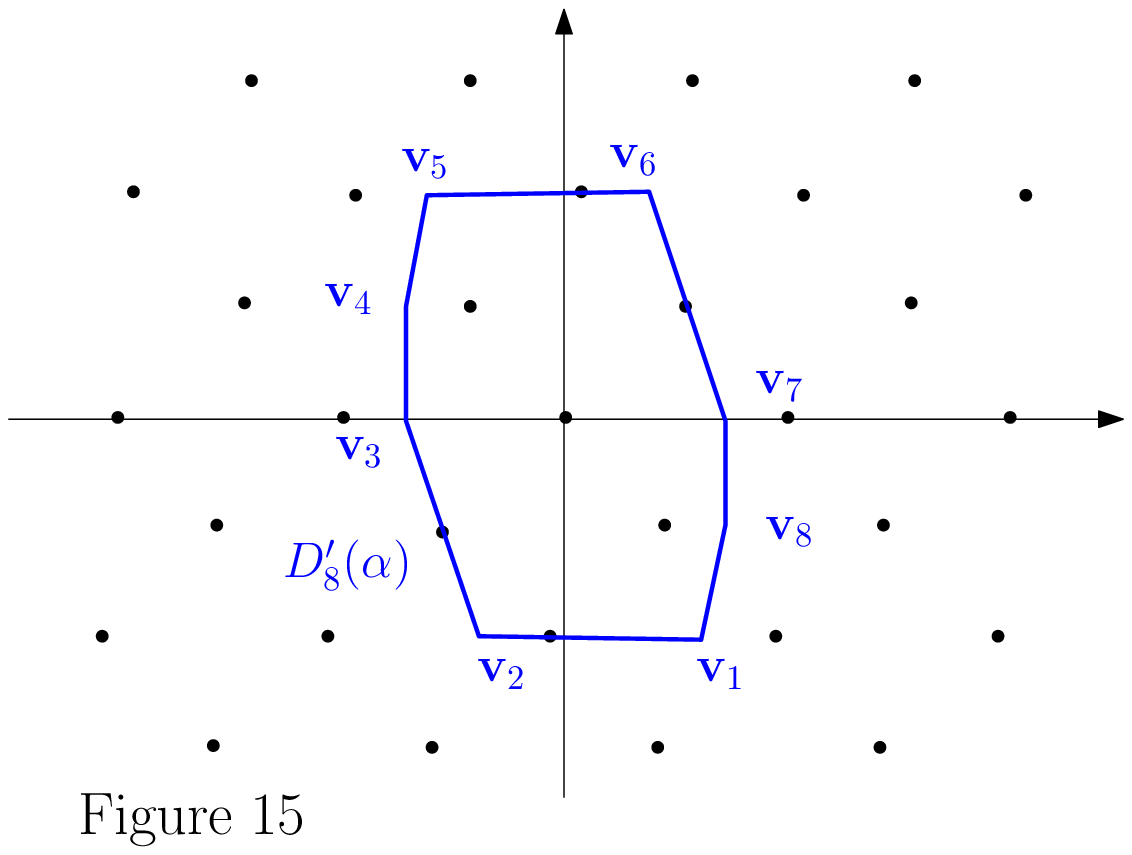}
\end{figure}

Similarly, let $D'_8(\beta )$ denote the octagon with vertices ${\bf v}_1=(2-\beta, -3),$ ${\bf v}_2=(-\beta, -3),$ ${\bf v}_3=(-2, -1),$ ${\bf v}_4=(-2, 1),$ ${\bf v}_5=-{\bf v}_1$, ${\bf v}_6=-{\bf v}_2$, ${\bf v}_7=-{\bf v}_3$ and ${\bf v}_8=-{\bf v}_4$, where $0<\beta \le 1$, and let $\Lambda (\beta )$ denote the lattice generated by ${\bf u}_1=(2,0)$ and ${\bf u}_2=(1+{\beta \over 2}, 2)$, as shown by Figure 16. By Lemma 2 it can be verified that $D'_8(\beta )+\Lambda (\beta )$ is indeed a five-fold translative tiling in $\mathbb{E}^2$.

\begin{figure}[!ht]
\includegraphics[scale=0.45]{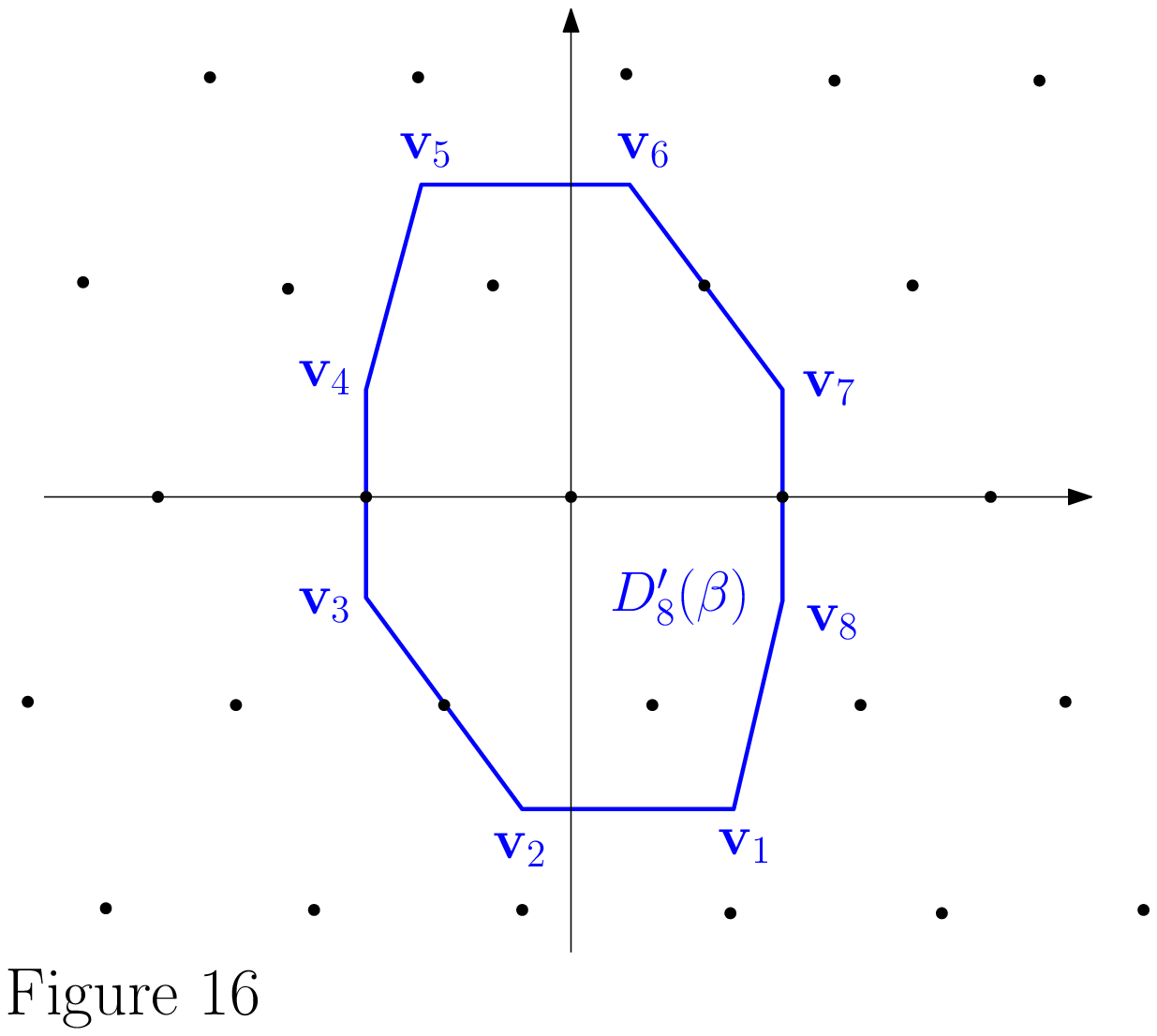}
\end{figure}

In fact, the two classes $D_8(\alpha )$ and $D'_8(\beta )$ shown in Figure 13 and Figure 16 respectively are equivalent under suitable linear transformations, as well as the two classes $D_8(\beta )$ and $D'_8(\alpha )$ shown by Figure 14 and Figure 15. Therefore, we have the following
theorem.

\medskip\noindent
{\bf Theorem 15 (Yang and Zong \cite{yz3}).} {\it A convex domain can form a five-fold translative tiling of the Euclidean plane if and only if it
can form a five-fold lattice tiling in $\mathbb{E}^2$.}

\vspace{0.6cm}\noindent
{\bf Acknowledgements.} The author is grateful to Prof. W. Zhang for calling his attention to Rao's announcement, to Profs M. Kolountzakis, S. Robins and G. M. Ziegler for helpful comments on related papers. This work is supported by 973 Program 2013CB834201.

\noindent
\bibliographystyle{amsplain}

\vspace{0.6cm}
\noindent
Chuanming Zong

\noindent
Center for Applied Mathematics

\noindent
Tianjin University, Tianjin 300072, China

\noindent
Email: cmzong@math.pku.edu.cn

\end{document}